\documentclass{elsart}
\usepackage{amssymb,amsmath,bbm}
\usepackage{graphicx}
\usepackage{epsfig}
\newtheorem{dl}{Theorem}

\newtheorem{yl}{Lemma}

\newtheorem{mt}{Proposition}
\newtheorem{dy}{Definition}
\newtheorem{rmk}{Remark}
\newtheorem{lz}{Example}

\newcommand{\R}{\mathbbmss{R}}
\newcommand{\n}{\mathcal{I}_n}

\newcommand{\yi}{\mathbf{1}}
\newcommand{\N}{\mathcal{N}}
\newcommand{\Z}{\mathbbmss{Z}_+}
\newcommand{\myA}{\mathcal{A}}
\newcommand{\mya}{\mathbf{a}}

\newcommand{\ssA}{{\mbox{\scriptsize\boldmath{$A$}}}}
\newcommand{\sA}{{\mbox{\boldmath{$A$}}}}
\newcommand{\sa}{{\mbox{\boldmath{$a$}}}}
\newcommand{\tul}{\check{\tau}_u}
\newcommand{\tuu}{\hat{\tau}_u}
\newcommand{\mml}{\check{m}}
\newcommand{\mmu}{\hat{m}}
\newcommand{\mmm}{\tilde{m}}
\newcommand{\tou}{\hat{\tau}_v}
\newcommand{\num}{\mbox{\boldmath{$n$}}}
\newcommand{\F}{\mathcal{F}}

\newcommand{\G}{\mathcal{G(A}(t))}
\newcommand{\Ge}{\mathcal{G}^0(t)}

\begin{document}

\begin{frontmatter}



\title{Asynchronous Consensus in Continuous-Time Multi-Agent Systems With Switching Topology and Time-Varying Delays \thanksref{funds}}


\author{Feng Xiao},
\ead{ fengxiao@pku.edu.cn}
\author{Long Wang}
\ead{longwang@pku.edu.cn}
\address{Intelligent Control Laboratory, Center for Systems and Control,
Department of Industrial Engineering and Management, and
Department of Mechanics and Space Technologies, College of
Engineering, Peking University, Beijing 100871, China}
\thanks[funds]{This work was supported by NSFC (60674050 and 60528007), National 973 Program (2002CB312200), and 11-5 project (A2120061303).}

\begin{abstract}
In this paper, we study asynchronous consensus problems of
continuous-time  multi-agent systems with discontinuous
information transmission. The proposed consensus control strategy
is implemented only based on the state information at some
discrete times of each agent's neighbors. The asynchronization
means that each agent's update times, at which the agent adjusts
its dynamics, are independent of others'. Furthermore, it is
assumed that the communication topology among agents is
time-dependent and the information transmission is with bounded
time-varying delays. If the union of the communication topology
across any time interval with some given length contains a
spanning tree, the consensus problem is shown to be solvable.  The
analysis tool developed in this paper is based on the nonnegative
matrix theory and graph theory. The main contribution of this
paper is to provide a valid distributed consensus algorithm that
overcomes the difficulties caused by unreliable communication
channels, such as intermittent information transmission, switching
communication topology, and time-varying communication delays, and
therefore has its obvious practical applications. Simulation
examples are provided to demonstrate the effectiveness of our
theoretical results.
\end{abstract}

\end{frontmatter}

\section{Introduction}
In recent years,  decentralized coordination of multi-agent
systems  has become an active area of research and attracted the
attention of multi-disciplinary researchers in a wide range
including system control theory, statistical physics, biology,
applied mathematics, and computer science. This is partly due to
its broad applications in cooperative control of unmanned aerial
vehicles, scheduling of automated highway systems, formation
control of satellite clusters, distributed optimization of
multiple mobile robotic systems, etc.

In  cooperative control of multiple agents, in order to accomplish
some complicated tasks or reach their common goals, groups of
dynamic agents in multi-agent/multi-robot systems   need to
interact with each other and eventually reach an agreement on
certain quantities of interest. Those problems are usually called
consensus problems, which are one of fundamental research topics
in decentralized control.

Consensus problems have been studied for a long time and the
formal investigation of them can be traced back to 1960's in the
field of management science and statistics (See DeGront
\cite{degroot} and references therein). In the field of systems
and control theory, the pioneering work was done by Borkar and
Varaiya\cite{V. Borkar and P. Varaiya} and Tsitsiklis and
Athans\cite{J. N. Tsitsiklis and M. Athans}, which is on
asynchronous consensus problems with an application in distributed
decision-making system. In \cite{T. Vicsek A. Czirok E. Ben Jacob
I. Cohen and O. Schochet}, Vicsek et al. proposed a simple but
interesting discrete-time model of multiple agents moving in the
plane. Each agent's motion is updated using a local rule based on
its own state and the states of its neighbors. The Vicsek model
can be viewed as a special case of a computer model mimicking
animal aggregation proposed in \cite{C. Reynolds} for the computer
animation industry. By using the graph theory and matrix theory,
Jadbabaie et al. provided a theoretical explanation of the
consensus property of the Vicsek model in \cite{A. Jadbabaie J.
Lin and A. S. Morse}, where each agent's set of neighbors changes
with time as the system evolves. The typical continuous-time
consensus model was presented by Olfati-Saber and Murray in
\cite{R. Olfati-Saber and R. M. Murray 1}, where the concepts of
solvability of consensus problems and consensus protocols were
first introduced. In \cite{R. Olfati-Saber and R. M. Murray 1},
under the assumptions that the dynamics of each agent is a scalar
continuous-time integrator and the communication among agents may
be unidirectional, Olfati-Saber and Murray used a directed graph
to model the communication topology
 and studied three agreement problems, namely, directed
networks with fixed topology, directed networks with switching
topology, and undirected networks with communication time-delays
and fixed topology. And it was assumed that the directed topology
is balanced and strongly connected. In \cite{W. Ren and R. W.
Beard}, Ren and Beard  extended the results of \cite{A. Jadbabaie
J. Lin and A. S. Morse} and \cite{R. Olfati-Saber and R. M.
Murray 1} and presented some improved conditions for state
agreement under dynamically changing directed interaction
topology, which is not necessarily balanced or strongly connected.
In the past several years, consensus problems of multi-agent
systems have been developing very fast and several research topics
have been addressed, such as agreement over random networks
\cite{Yuko,A. V. Savkin}, asynchronous information consensus
\cite{L. Fang,L. Fang 2}, dynamic consensus \cite{D. P. Spanos
dynamic consensus}, consensus filters \cite{Reza Olfati-Saber2},
networks with general communication structures \cite{xfphysicaA},
networks with nonlinear agreement protocols \cite{L. Moreau 1},
and networks with switching topology and time-delays \cite{A.
Jadbabaie J. Lin and A. S. Morse,R. Olfati-Saber and R. M. Murray
1,W. Ren and R. W. Beard,L. Moreau 1,Dongjun Lee,L. Moreau 2,feng
xiao,Tanner}. For more details, see the survey \cite{W. Ren
survey} and paper \cite{R. Olfati-Saber proceeding}, where a
theoretical frame-work for analysis of consensus algorithms for
multi-agent networked systems was provided. In addition, it is
necessary  to emphasize that flocking of agents and swarms\cite{R.
Olfati-Saber flocking,H. Shi pD,Wang,T. Chu} and formation control
of vehicles\cite{zhiyun lin1,zhiyun lin2,lafferriere} are two
active areas where many useful results obtained in  consensus
problems have been successfully applied. And in \cite{R.
Olfati-Saber flocking}, the authors provided the first proof of
convergence of Reynolds' rules on the basis of the convergence of
consensus algorithms  in \cite{R. Olfati-Saber and R. M.  Murray
1}.

In this paper, we propose a distributed asynchronous consensus
control strategy that is only based on the state information of
each agent's neighbors at some discrete times. This is partly
motivated by the work of Olfati-Saber and Murray \cite{R.
Olfati-Saber and R. M. Murray 1} and the difficulties encountered
in the implementation of the typical continuous-time protocols
proposed in \cite{R. Olfati-Saber and R. M.  Murray 1}. In the
realistic networks of agents, we may face the following problems:
\begin{enumerate}
    \item Communication topology is always changing;
    \item The received information is often with time-delays, and
    furthermore, the delays may be (randomly) time-varying and unknown;
    \item Due to long distance transmission, unreliable
information channels,  and limited bandwidth of networks, the
continuity of state information of each agent's
    neighbors cannot be ensured.
\end{enumerate}

Here, we summarize some important  closely related works. First,
in continuous-time systems, none of the existing results guarantee
the stability of the consensus protocols proposed by Olfati-Saber
and Murray, in the presence of switching topology and time-varying
delays. It was often assumed that the available topologies are
finite\cite{A. Jadbabaie J. Lin and A. S. Morse,R. Olfati-Saber
and R. M. Murray 1,W. Ren and R. W. Beard} and the delays are
constant\cite{R. Olfati-Saber and R. M. Murray 1,L. Moreau
1,Dongjun Lee}. In fact, the stability analysis of the protocols
in \cite{R. Olfati-Saber and R. M. Murray 1} in the presence of
time-varying delays is challenging and even impossible in theory,
and is unnecessary in applications.

In the study of discrete-time systems, Tanner and Christodoulakis
in \cite{Tanner} studied a discrete-time model with fixed
undirected topology and assumed that all agents transmit their
state information in turn. Consequently, outdated information may
be used and the equivalent augmented system becomes a periodically
switched system, which can be viewed as a multi-agent system with
switching topology. In \cite{L. Fang}, Fang and Antsaklis studied
the case with switching topology and time-dependent delays by an
asynchronous system with fixed topology. However since the
possible topologies are generated by the fixed topology of the
asynchronous system, the ``switching" topology is not really
switching. By matrix theory, Xiao and Wang derived some sufficient
conditions for the solvability of consensus problems of
discrete-time systems with switching topology and time-varying
delays in \cite{feng xiao}, but it was assumed that all available
topologies are finite.

In addition, all the proposed continuous-time consensus algorithms
 depend on continuous state signals. In engineering applications,
continuous signals  require large bandwidth of networks, and
furthermore, in many cases, are not  available. Therefore, we need
to devise some continuous-time consensus algorithms  that do not
continuously depend on the external state information.

Our proposed consensus control strategy is built upon very general
assumptions. The communication topology is switching among an
infinite set of weighted directed graphs, communication delays are
time-varying, and state information transmission is allowed to be
intermittent. Moreover, our control strategy is also an
asynchronous one, which means that each agent's update actions are
independent of others'. Each agent adjusts its dynamics
independently, which is inherent in distributed control systems.
It is important to mention that asynchronous consensus problems
were also studied by \cite{ V. Borkar and P. Varaiya,J. N.
Tsitsiklis and M. Athans} and \cite{L. Fang}, where several
asynchronous consensus algorithms were given. By using nonnegative
matrix theory and graph theory, especially the properties of
scrambling matrices, we provide some sufficient (and necessary)
conditions for the convergence of our consensus control strategy.

This paper is organized as follows. Section II presents some basic
definitions and results in matrix theory and graph theory. Section
III  formulates the problem to be studied. Convergence analysis
and the technical proof are performed in Sections IV and V,
respectively. In Section VI, simulation examples are presented.
Finally, concluding remarks are stated in Section VII.

\section{Preliminaries}

This section  presents some definitions and results in matrix
theory and graph theory that will be used in this paper \cite{R.
Horn and C. R. Johnson,C. Godsil and G. Royal}.

Let $\n=\{1, 2, \cdots, n\}$, $\Z$ be the set of nonnegative
integers, and $\yi=[1,1, \cdots, 1]^T$ with compatible dimensions.
Given $A=[a_{ij}] \in \R^{n\times r}$,   $A$ is said to be {\it
nonnegative},  $A\geq 0$, if all its entries $a_{ij}$ are
nonnegative.  $A$ is said to be {\it positive}, $A>0$,  if all its
entries $a_{ij}$ are positive. Let $B\in \R^{n\times r}$. We write
$A\geq B$ if $A-B\geq 0$, and $A> B$ if $A-B> 0$.  A { nonnegative
square matrix} $A$ with the property that all its row sums are
$+1$ is said to be a {\it stochastic matrix}. Throughout this
paper, we let $\prod_{i=1}^k A_i=A_kA_{k-1}\cdots A_1$ denote the
left product of matrices. A $n\times n$ stochastic matrix $A$ is
called indecomposable and aperiodic (SIA) (or {\it ergodic}) if
there exists $f\in\R^n$ such that $\lim_{k\to \infty}A^k=\yi f^T$.

Directed graphs  will be used to model the communication
topologies among agents. A {\it directed graph} $\mathcal{G}$
consists of a vertex set $\mathcal{V}(\mathcal{G})=\{ v_1, v_2,
\cdots, v_n\}$ and an edge set $\mathcal{E}(\mathcal{G})\subset\{
(v_i, v_j): v_i, v_j \in \mathcal{V}(\mathcal{G})\}$, where an
edge is an ordered pair of vertices in $\mathcal{V}(\mathcal{G})$
(Here, we allow for self-loops, namely, the edges with the same
vertices). The set of {\it neighbors} of vertex $v_i$ in
$\mathcal{G}$ is denoted by $\mathcal{N(G,}v_i)=\{ v_j: (v_j,v_i)
\in \mathcal{E}(\mathcal{G}), j\not=i\}$. The associated index set
of the neighbors is denoted by $\mathcal{N(G},i)=\{j:v_j\in
\mathcal{N}(\mathcal{G}, v_i)\}$. If $(v_i,v_j)$ is an edge of
$\mathcal{G}$, $v_i$ and $v_j$ are defined as the parent and child
vertices, respectively. A {\it subgraph} $\mathcal{G}_s$ of a
directed graph $\mathcal{G}$ is a directed graph such that the
vertex set $\mathcal{V}(\mathcal{G}_s) \subset
\mathcal{V}(\mathcal{G})$ and the edge set
$\mathcal{E}(\mathcal{G}_s) \subset \mathcal{E}(\mathcal{G})$.  If
$\mathcal{V}(\mathcal{G}_s) = \mathcal{V}(\mathcal{G})$, we call
$\mathcal{G}_s$ a {\it spanning subgraph} of $\mathcal{G}$. For
any $v_i, v_j\in\mathcal{V}(\mathcal{G}_s)$, if $(v_i,v_j)\in
\mathcal{E}(\mathcal{G}_s)$ if and only if $(v_i,v_j)\in
\mathcal{E}(\mathcal{G})$, $\mathcal{G}_s$ is called an {\it
induced subgraph}.  In this case, $\mathcal{G}_s$ is also said to
be induced by $\mathcal{V}(\mathcal{G}_s)$. A {\it path} in a
directed graph $\mathcal{G}$ is a sequence $v_{i_1}, \cdots,
v_{i_k}$ of vertices such that $(v_{i_s}, v_{i_{s+1}})\in
\mathcal{V(G)}$ for $s=1, \cdots, k-1$.   A directed graph
$\mathcal{G}$ is {\it strongly connected} if between every pair of
distinct vertices $v_i, v_j$ in $\mathcal{G}$, there is a path
that begins at $v_i$ and ends at $v_j$ (that is, from $v_i$ to
$v_j$).  A {\it directed tree} is a directed graph, where every
vertex, except one special vertex without any parent, which is
called the {\it root vertex}, has exactly one parent, and the root
vertex can be connected to any other vertices through paths. A
{\it spanning tree} of $\mathcal{G}$ is a directed tree that is a
spanning subgraph of $\mathcal{G}$. We say that a graph has (or
contains) a spanning tree if a subset of the edges forms a
spanning tree.
 A {\it weighted
directed graph} $\mathcal{G}(A)$ is a directed graph $\mathcal{G}$
plus a nonnegative {\it weight matrix} $A=[a_{ij}]\in \R^{n\times
n}$ such that $(v_i, v_j)\in \mathcal{E(G)} \iff a_{ji}>0$. And
$a_{ji}$ is called the {\it weight} of edge $(v_i, v_j)$.

\section{Problem Formulation}
We suppose that the system studied in this paper consists  of $n$
autonomous agents, e.g., birds, robots, etc.,  labeled $1$ through
$n$. All these agents share a common state space $\R$. Each agent
adjusts its current state based upon the information received from
other agents that are defined as { neighbors} of this agent. We
use a weighted directed graph $\mathcal{G(A}(t))$ to represent the
communication topology or information flow, where
$\mathcal{A}(t)=[\mya_{ij}(t)]\in\R^{n\times n}$ is a nonnegative
matrix.  The appearance of parameter $t$ implies that the
communication topology may be dynamically changing. Agent $i$ is
represented by vertex $v_i$. Edge $(v_j, v_i)\in\mathcal{E}(\G)$
corresponds an available information channel from agent $j$ to
agent $i$. If agent $i$ receives information from agent $j$ at
time $t$, then there exists an edge from vertex $v_j$ to vertex
$v_i$, i.e., $(v_j,v_i)\in\mathcal{E}(\G)$. And the {\it
neighbors} of agent $i$ are those agents whose information is
received by agent $i$ at time $t$. The associated index set of the
neighbors is denoted by $\N(t,i)$. Notice that because of the
existence of communication time-delays, the index set $\N(\G,i)$
may not be equal to $\N(t,i)$. We will discuss them latter.

Let $x_i\in\R$ denote the state of agent $i$ and let $x=[x_1, x_2,
\cdots, x_n]^T$. Then the whole system can be generally
represented by the continuous-time model $\dot{x}(t)=f(t,u(t))$ or
by the discrete-time model $x(t+1)=f(t,u(t))$, where $u(t)$ is a
state feedback. If for any initial state, $x(t)$ converges to some
equilibrium point $x^*$ (dependent on the initial state) such that
$x^*_i=x^*_j$ for all $i, j\in \n$, as $t\to \infty$, then we say
that this system {\it solves a consensus problem}\cite{R.
Olfati-Saber and R. M.  Murray 1} (or {\it has the consensus
property}). Let $\chi:\R^n\to \R$ be a function of $n $ variables
$x_1, x_2, \cdots, x_n$.  For any initial state
$x(0)$,\footnote{In this paper, consensus functions are only
related to the systems without time-delays.} if $x^*=\yi
\chi(x(0))$, then we say that this system {\it solves the
$\chi$-consensus problem}, and the function $\chi$ is called the
{\it consensus function}. The common value of $x^*_i$ is called
the {\it group decision value}.

\subsection{The Model}

For agent $i$, $i\in\n$, we assume that it receives or detects its
neighbors' states at {\it update times} $t_0^i, t_1^i,\cdots,
t_k^i, \cdots$, which can be seen as a real number sequence and
are denoted by $\{t_k^i\}$. We assume that $\{t_k^i\}$ satisfies
the following
assumptions:\\

    (A1) For any $k\in\Z$, $0<\tul\leq t_{k+1}^i-t_{k}^i\leq
    \tuu,$
      where $\tul, \tuu\in\R$;

    (A2) $t_0^i=0$.\\

The simple reason for the calling of ``update times" is that the
neighbors' information known by  agent $i$ or the dynamics of
agent $i$ is updated at those times. The existence of lower bound
$\tul$ of time intervals between any two consecutive update times
is just to guarantee the validity of our consensus protocols
\eqref{xt1} and \eqref{xtdelay}. If there does not exist lower
bound, it will be hard to analyze protocols \eqref{xt1} and
\eqref{xtdelay} in theory and it is also unnecessary in
applications. If there does not exist an upper bound $\tuu$ of
$t_{k+1}^i-t_k^i$, it will be difficult for the states of agents
to reach consensus (See Example \ref{slz2}). For Assumption (A2),
we make it solely for the convenience of our theoretical analysis;
otherwise, the main results of our paper is still obtainable.

 If agent $i$
receives  the state information of its neighbors at $t_k^i$, then
agent $i$ is assumed to take the following dynamics in  time
interval $[t_k^i, t_{k+1}^i)$
\begin{equation}\label{xt1}
    \dot{x}_i(t)=\left\{
                   \begin{array}{ll}
                     \sum_{j\in\N(t_k^i,i)}\frac{\mya_{ij}(t_k^i)}{\sum_{j\in\N(t_k^i, i)}\mya_{ij}(t_k^i)}(x_j(t_k^i)-x_i(t)),& \mbox{if\ } \N(t_k^i,i)\not=\phi; \\
                    0, & \mbox{otherwise.} \\
                   \end{array}
                 \right.
\end{equation}

\begin{rmk}
 The continuous-time model studied in \cite{R.
Olfati-Saber and R. M. Murray 1} is that
$\dot{x}_i(t)=\sum_{j\in\N(t,i)}\mya_{ij}(t)$ $ (x_j(t)-x_i(t))$.
In our discrete-communication framework, this system will turn
into $\dot{x}_i(t)=\sum_{j\in\N(t,i)}\mya_{ij}$
$(t_k^i)(x_j(t_k^i)-x_i(t))$ accordingly, where $t\in [t_k^i,
t_{k+1}^i)$. Seemingly, our model \eqref{xt1} is a special case of
this one. But if system \eqref{xt1} solves a consensus problem,
then the afore-mentioned system will also have the consensus
property. That will be clear in our analysis process.
\end{rmk}

It is well known that communication time-delays exist extensively
in networks. Therefore, it is reasonable to assume that there
exist communication time-delays in the information transmission.
For agent $i$, $i\in\n$, we assume that the information received
by agent $i$ from agent $j$ is with time-delay $\tau^k_{ij}$ at
update time $t_k^i$, and then the considered system turns into
\begin{equation}\label{xtdelay}
    \dot{x}_i(t)=\left\{
                   \begin{array}{ll}
                     \sum_{j\in\N(t_k^i,i)}\frac{\mya_{ij}(t_k^i)}{\sum_{j\in\N(t_k^i, i)}\mya_{ij}(t_k^i)}(x_j(t_k^i-\tau^k_{ij})-x_i(t)),& \mbox{if\ } \N(t_k^i,i)\not=\phi; \\
                    0, & \mbox{otherwise.} \\
                   \end{array}
                 \right.
\end{equation}
where $t\in[t_k^i, t_{k+1}^i)$.

We make the following assumption about system \eqref{xtdelay}
together with
Assumption (A1, A2):\\

(A3) $0\leq\tau^k_{ij}\leq\tau_d=K\tul$, where $i\in\n$,
$j\in\N(t_k^i,i), k,K\in\Z$;\\

\begin{rmk}
Consider system \eqref{xtdelay}. Generally, the time-delays may be
unknown for each each agent when they receive information. If the
data is time stamped, the delays will be detectable. In addition,
if communication time-delays only satisfy Assumption (A3), for any
$i\in \n$, $j\in\N(t_k^i,i)$, the state information of agent $j$,
$x_j(t_{k}^i-\tau_{ij}^{k})$, received by agent $i$ at update time
$t_{k}^i$ may be outdated compared
 with state information $x_j(t_{k'}^i-\tau_{ij}^{k'})$  received previously, i.e.,
 $t_{k}^i-\tau_{ij}^{k}<t_{k'}^i-\tau_{ij}^{k'}$, where $k'<k$.  If the delays are detectable, we can suppose
 that agent $i$ always uses the most recent data of its neighbors,
 that is, if there exists $k'<k$ such that $t_{k}^i-\tau_{ij}^{k}<t_{k'}^i-\tau_{ij}^{k'}$, then
 agent $i$ replaces the state information  $x_j(t_{k}^i-\tau_{ij}^{k})$ by
 $x_j(t_{k^*}^i-\tau_{ij}^{k^*})$ in time interval $[t_{k}^i,
 t_{k+1}^i)$, where $k^*=\mbox{arg}\max_{k'<k}(t_{k'}^i-\tau_{ij}^{k'})$. We call this control strategy {\it the-most-recent-data
 strategy}. This control strategy can get better convergence rate (See Example
 \ref{lzdelay}).
\end{rmk}

\begin{figure}[thpb]
\centering
      \includegraphics[scale=1]{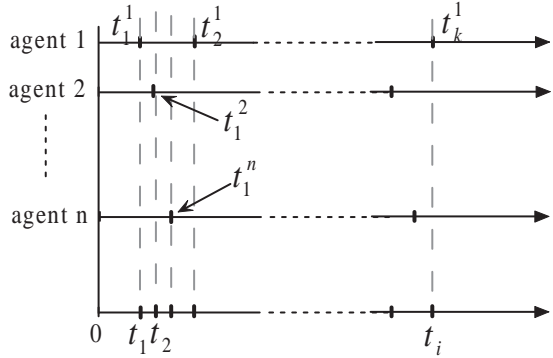}
           \caption{Update times in asynchronous systems without time-delays.}
      \label{fig2}
\end{figure}

According to different properties of update times, we classify
systems \eqref{xt1} and \eqref{xtdelay} as synchronous or as
asynchronous systems.
\begin{dy}[Synchronous and asynchronous systems]
We say that system \eqref{xt1} (or \eqref{xtdelay}) is {\it
synchronous} if for any $i,j\in\n$, $\{t_k^i\}=\{t_k^j\}$, i.e.,
for any  $k\in\Z$, $t_k^i=t_k^j$. We say that system \eqref{xt1}
(or \eqref{xtdelay}) is {\it asynchronous} if for any $i,j\in\n$,
$i\not=j$, $\{t_k^i\}$ is independent of $\{t_k^j\}$, i.e., agents
may not adjust their dynamics at the same time (See Fig.
\ref{fig2}).
\end{dy}

In realistic networks, it is difficult for all agents to be
synchronous on update actions and therefore we mainly discuss the
asynchronous consensus property of systems \eqref{xt1} and
\eqref{xtdelay}.

\subsection{Communication Topology}

Since our consensus protocols \eqref{xt1} and \eqref{xtdelay} only
depend on discrete state information, we are not concerned with
the actual communication topology $\G$ outside those update times.
And we give the following definition of $\Ge$, which is different
from the actual communication topology and is also called
communication topology.

\begin{dy}[Communication topology]\label{dytopology} Suppose that $\Ge$ is with the same vertex set as $\G$.
 For any $i\in\n$, $k\in\Z$,
if agent $i$ receives the state information of agent $j$ at time
$t_k^i$, $j\not=i,j\in\n$, then $(v_j,v_i)\in\Ge$, and if not,
$(v_j,v_i)\not\in\Ge$, where $t\in [t_k^i, t_{k+1}^i)$. In
addition, there are no self-loops in $\Ge$.
\end{dy}

In Definition \ref{dytopology}, $\Ge$ is only compliant with the
actual topology $\G$ on those edges corresponding to
$\mya_{ij}(t_k^i)$, $i\in\n, j\in \N(t_k^i,i)$. This way can
facilitate our theoretical analysis and does not affect the
dynamics of agents.

\begin{mt}
$\N(t_k^i,i)=\N(\mathcal{G}^0(t_k^i),i)\subset\N(\mathcal{G(A}(t_k^i)),i)$.
\end{mt}

The above proposition follows from the fact that the available
channel $(v_j,v_i)$ at $t_k^i$ can not ensure that agent $i$ can
receive the state information of agent $j$ at time $t_k^i$  if
there exists time-delay, and  agent $i$ cannot receive any
information from agent $j$   if $(v_j,v_i)$ does not exist.

\begin{rmk}
The reason for the introduction of $\Ge$ is that $\Ge$ is more
important than $\G$ for our control strategies. We can see that
$\Ge$ is directly connected with the successful information
transmission in the networks. $\mathcal{E}(\G)$ represents all
available communication channels, while $\mathcal{E}(\Ge)$ stands
for the communication channels through which the information has
been successfully received.
\end{rmk}

Here, for simplicity of presentation, we introduce  another matrix
$\sA(t)=[\sa_{ij}(t)]$, where $\sa_{ij}(t)=\left\{
                      \begin{array}{cc}
                        \frac{\mya_{ij}(t_k^i)}{\sum_{j\in\N(t_k^i,i)}\mya_{ij}(t_k^i)}, &j\not=i \\
                        0, & j=i \\
                      \end{array}
                    \right.
$ if $\N(t_k^i,i)\not=\phi$,  and $\sa_{ij}(t)=\left\{
                      \begin{array}{cc}
                        0, &j\not=i \\
                        1,& j=i \\
                      \end{array}
                    \right.
$ if $\N(t_k^i,i)=\phi$, $t\in[t_k^i,t_{k+1}^i)$. Obviously,
$\sA(t)$ is stochastic.

\begin{mt}\label{mtgraphG}
For any $t\geq 0$, $\N(\mathcal{G}(\sA(t)),i)=\N(\Ge,i)$. And if
we ignore the weight of each edge in $\mathcal{G}(\sA(t))$ and
self-loops in $\mathcal{G}(\sA(t))$, then $\mathcal{G}(\sA(t))$
and $\Ge$ represent the same graph.
\end{mt}

Since $\G$ may be dynamically  changing, we should investigate
all possible directed graphs. Because of the finite number of
vertices, there are at most $2^{n \times n}$ different kinds of
directed graphs. Let $\Gamma_{\mathcal{G}}$ denote the set of all
those directed graphs. Assume the existence of real number
$\mathbf{\hat{a}}\geq \mathbf{\check{a}}>0$, such that
$\mathbf{\check{a}}\leq\mya_{ij}(t)\leq\mathbf{\hat{a}}$ if
$\mya_{ij}(t)\not=0$. Consequently all possible $\sA(t)$
constitute a compact set\footnote{The set of all $r\times s$
matrices can be viewed as the metric space $\R^{rs}$ and compact
sets are equivalent to bounded closed sets.}, denoted by
$\Gamma_{\ssA}$, and if $\sa_{ij}(t)\not=0$, then
$\frac{\mathbf{\check{a}}}{(n-1)\mathbf{\hat{a}}}< \sa_{ij}(t)\leq
1$. It is important to note that all possible communication
topologies are infinite if we take the weight of each edge into
account.

Finally, we present the notion of the union of graphs that will be
used in the remainder of this paper.
\begin{dy}[Union of graphs]
The union of graph $\Ge$ across time interval $[t^0, t^0+T]$ is a
directed graph with the same vertex set as $\Ge$ and the edge set
$\bigcup_{t'\in[t^0,t^0+T]}\mathcal{E}(\mathcal{G}^0(t'))$. The
union of graph $\Ge$ on the time set $\{s_1, s_2, \cdots, s_k\}$
is a directed graph with the same vertex set as $\Ge$ and the edge
set $\bigcup_{i=1}^k\mathcal{E}(\mathcal{G}^0(s_i))$.
\end{dy}

\section{Convergence Results}
This section presents the main result of this paper. As a
preparation for the study of the general case, we first
investigate two relatively simple cases: the synchronous system
with fixed topology in the absence of time-delays and the
asynchronous system without time-delays. Our  approach is to
transpose the continuous-time systems into their  discrete-time
counterparts, which possess the same consensus property as the
original systems. The obtained discrete-time systems have some
special structures and the proof of their consensus property is
postponed to the next section.

 Some notations are used in this section. If the studied systems
are free of time-delays, we let $\{t_k\}$ $=$
$\{t_k^i,i\in\n,k\in\Z\}$, and if the systems are with
time-delays, let $\{t_k\}=\{t_k^{i} \mbox{\ or\ }
t_k^i-\tau_{ij}^k, i\in\n, k\in\Z, j\in \N(t_k^i,i)\}$, such that
$t_0=0$ and $t_{k+1}>t_k$. Let $\tau_k=t_{k+1}-t_k$, $k\in\Z$, and
let $\Lambda(A)=\{B=[b_{ij}]\in\R^{n\times n}: \mbox{for any }
i,j\in\n, b_{ij}=a_{ij},$ $\mbox{or\ } b_{ij}=0\}$, where
$A=[a_{ij}]\in \R^{n\times n}$. Obviously $\Lambda(A)$ is a finite
set. Let $\Pi(m,t)$ denote the set of matrices
\begin{equation}\label{mat}
\left[
                            \begin{array}{ccccc}
                              e^{-h}I+(1-e^{-h})A_1(t) & (1-e^{-h})A_2(t) & \cdots & (1-e^{-h})A_{m-1}(t) & (1-e^{-h})A_m(t) \\
                              I & 0 & \cdots & 0 & 0 \\
                               0 & I& \cdots & 0 & 0 \\
                              \vdots & \vdots & \ddots & \vdots & \vdots \\
                              0 & 0 & \cdots & I &0 \\
                            \end{array}
                          \right]_{mn\times mn},
\end{equation}
where $0<h\leq \tuu$, $ A_1(t),\cdots,A_m(t)\in\Lambda(\sA(t))$,
and $ A_1(t)+\cdots+A_m(t)=\sA(t)$. We denote the above matrix
\eqref{mat} by $\pi(h, A_1(t),\cdots,A_m(t))$.

\subsection{The Synchronous Consensus}
We first restrict our attention to the simplest case and we study
the synchronous system with time-invariant topology in the absence
of time-delays. Since there exist no communication delays, we
assume that for any $i\in\n$,
$j\in\mathcal{N}(\mathcal{G(A}(t_k^i)),i)$, agent $i$ can receive
state information from agent $j$  at time $t_k^i$. Therefore
$\mathcal{N}(\mathcal{G(A}(t)),i)=\mathcal{N}(\mathcal{G}(\sA(t)),i)$,
$\N(t_k^i,i)=\N(\mathcal{G(A}(t_k^i)),i) $, and $\sA(t)$ is
time-invariant. For simplicity, notations $\myA$ and $\sA$ are
used instead of $\myA(t)$ and $\sA(t)$. Apparently,
$\{t_k^i\}=\{t_k\}$. Rewriting system \eqref{xt1} yields that for
any $i\in\n$,  $t\in [t_k, t_{k+1})$
\begin{equation}\label{sxt1}
\dot{x}_i(t)=\left\{
  \begin{array}{ll}
    -x_i(t)+\sum_{j\in\N(\mathcal{G(A}), i)}\sa_{ij}x_j(t_k), & \mbox{if\ } \N(\mathcal{G(A)},i)\not=\phi;\\
     0,& \mbox{otherwise.}\\
  \end{array}
\right..
\end{equation}
Then we have
\begin{equation}\label{sxt2}
    \begin{split}
    x(t_k+h)&=(e^{-h}I+(1-e^{-h})\sA)x(t_k)\\
            &=((1-e^{-h})(\sA-I)+I)x(t_k),
    \end{split}
\end{equation}
where $0\leq h \leq\tau_k$ and $I$ is the identity matrix with
compatible dimensions.

The following theorem characterizes the consensus property of
synchronous system \eqref{sxt1}.
\begin{dl}\label{sdl}
System \eqref{sxt1} solves a consensus problem if and only if
$\mathcal{G(A)}$ has a spanning tree. In addition, the group
decision value is uniquely determined by $\mathcal{G(A)}$ and the
initial state.
\end{dl}

The sufficiency of the first part is a direct consequence of
Theorem \ref{dlmain}, and we only prove the necessity and the
second statement.

Proof: The necessity is shown first. If $\mathcal{G(A)}$ has not
any spanning tree, then there will be several subsystems, among
which there will not be information communication, and thus system
\eqref{sxt1} will  not solve any consensus problem. Next, we prove
the second statement. Because $\sA$ is a stochastic matrix, there
exists an $f\in\R^n$, $f\geq 0$, such that $f^T\yi=1$ and
$f^T\sA=f^T$( Theorem 1, \cite{xfjmaa}). Let $\chi(x)=f^Tx$.
\[\chi(x(t_k+h))=f^T((1-e^{-h})(A-I)+I)x(t_k)=f^Tx(t_k)=\chi(x(t_k)).\]
Therefore, $\chi(x)$ is time-invariant.

Suppose that the final state is $\yi a$, $a\in\R$. We have that
$\lim_{t\to\infty}x(t)=\yi a$ and
$\lim_{t\to\infty}\chi(x(t))=\chi(\yi a)= f^T\yi a =a$. It follows
that the group decision value $a= \chi(x(0))$, which is uniquely
determined by $\mathcal{G(A)}$ and the initial state, and the
consensus function is $\chi(x)$.

Now, we  give an example to show that the assumption of the
existence of upper bound $\tuu$ of update intervals  is necessary.
\begin{lz}[{Counterexample}]\label{slz2}
Consider the synchronous case and let $n=2$ and
$\mathcal{A}=\left[
                                                 \begin{array}{cc}
                                                   0 & 1 \\
                                                   1 & 0 \\
                                                 \end{array}
                                               \right].
$ If update intervals $\tau_0=\ln 8$, $\tau_1=\ln 16$, $\cdots$,
$\tau_{k}=\ln 2^{k+3}$, $\cdots$, and the initial state $x(0)=[1,
-1]^T$, then we obtain that
\[\lim_{k\to\infty}|t_{k+1}-t_k|=\lim_{k\to\infty}\ln2^{k+3}=\infty,\]
and \begin{equation*}
    \begin{split}
    x_{1}(t_0)-x_2(t_0)&=2;\\
    x_2(t_1)-x_1(t_1)&=2(1-\frac{1}{4});\\
    x_1(t_2)-x_2(t_2)&=2(1-\frac{1}{4})(1-\frac{1}{8});\\
    \vdots
    \end{split}
\end{equation*}
For any $k\geq 1$, we have
\begin{equation*}
    \begin{split}
|x_1(t_k)-x_2(t_k)|&=2(1-\frac{1}{4})(1-\frac{1}{8})\cdots(1-\frac{1}{2^{k+1}})\\
&\geq 2(1-\frac{1}{4}-\frac{1}{8}-\cdots-\frac{1}{2^{k+1}}).
    \end{split}
\end{equation*}
Therefore, $\lim_{k\to\infty}|x_1(t_k)-x_2(t_k)| \geq 2(1-0.5)=1$.
Consequently, the states of agent $1$ and agent $2$ will never
reach consensus.
\end{lz}

\subsection{The Asynchronous Consensus}

\subsubsection{The Case Without Time-Delays}

Consider the asynchronous case of system \eqref{xt1}. Lemma
\ref{ylwithoutdelay} gives a proper characterization of update
times, which is important in the forthcoming analysis. From
Assumption (A1), we have the following fact
\begin{yl}\label{ylwithoutdelay}
Consider system \eqref{xt1}. For any $i\in\n$, the number of
elements in set $\{t_j:t_j\in[t_k^i,t_{k+1}^i)\}$ is no greater
than $(\lfloor \frac{\tuu}{\tul} \rfloor+1)(n-1)+1$ denoted by
$\mml$, where $\lfloor\frac{\tuu}{\tul}\rfloor$ is the maximum
integer no greater than $\frac{\tuu}{\tul}$.
\end{yl}

Proof: By Assumption (A1), we have that $t^i_{k+1}-t^i_k\leq
\tuu$. For any $j\in\n, j\not=i$, agent $j$ updates state
information of its neighbors  at most
$\lfloor\frac{\tuu}{\tul}\rfloor+1$ times in time interval
$(t^i_k, t_{k+1}^i)$. And there are $n-1$ possible $j$s. Take
$t_k^i$ into account and therefore the number of elements in
$\{t_j:t_j\in[t^i_k, t_{k+1}^i)\}$ is not greater than $(\lfloor
\frac{\tuu}{\tul} \rfloor+1)(n-1)+1$.

For agent $i$, $i\in\n$, and update time $t_k$, $k\geq \mml-1$,
there exists $s\in\Z$ such that $t_s^i\leq t_k< t_{k+1}\leq
t_{s+1}^i$. By solving equation equation \eqref{xt1}, we have that
\[x_i(t_{k+1})=e^{-\tau_k}x_i(t_{k})+(1-e^{-\tau_k})\sum_{j\in\N(t_s^i,i)\cup\{i\}}\sa_{ij}(t_s^i)x_{j}(t_s^i),\]
where $\sa_{ij}(t_s^i)=\sa_{ij}(t_k)$.

 By Lemma \ref{ylwithoutdelay}, $t_s^i\geq t_{k-\mml+1}$.
Let $y(k)=[x(t_k)^T,x(t_{k-1})^T,\cdots,x(t_{k-\mml+1})^T]^T$,
where $k\geq \mml-1$. From the above discussion, there exists a
matrix $\pi(\tau_k,\sA_{1}(t_k),\cdots,$
$\sA_{{\mml}}(t_k))\in\Pi(\mml,t_k)$ such that
\begin{equation}\label{xtdescrete2}
    y(k+1)=\pi(\tau_k,\sA_{1}(t_k),\cdots,\sA_{\mml}(t_k))y(k).
\end{equation}

\begin{mt}\label{mt1}
System \eqref{xt1} solves a consensus problem if and only if
system \eqref{xtdescrete2} solves a consensus problem.
\end{mt}
Proof: The necessity follows from the definition of state variable
$y(k)$. Assume that  system \eqref{xtdescrete2} solves a consensus
problem. Let $a\in\R$ such that $\lim_{k\to\infty}y(k)=\yi a$. For
any $i\in\n$,
 $x_i(t)-x_i(t_k^i)=(1-e^{-(t-t_k^i)})(\sum_{j=1}^n\sa_{ij}(t_k^i)x_j(t^i_k)-x_i(t^i_k))$,
 where $t_k^i < t \leq t_{k+1}^i$. Since $t^i_k\to\infty$ as
 $t\to\infty$, $\lim_{t\to\infty}(\sum_{j=1}^n\sa_{ij}(t_k^i)x_j(t_k^i)-x_i(t_k^i))=\lim_{t\to\infty}
 \sum_{j=1}^n\sa_{ij}(t_k^i)(x_j(t_k^i)-x_i(t_k^i))=0$ (for $\lim_{t\to\infty}(x_j(t_k^i)-x_i(t_k^i))=a-a=0$),
   that is,  $\lim_{t\to\infty}(x_i(t)-x_i(t_k^i))=0$.
 Therefore, for any $i\in \n$, $\lim_{t\to\infty}x_i(t)=a$, and
 system \eqref{xt1} solves a consensus problem.

The following theorem is also a consequence of Theorem
\ref{dlmain}.
\begin{dl}
If there exists $T\geq0$ such that for all $t^0\geq 0$, the union
of graph $\Ge$ across interval $[t^0, t^0+T]$ contains a spanning
tree, then
 system \eqref{xt1} solves a consensus problem.
\end{dl}

\subsubsection{The Case With Time-Delays}

As a parallel with Lemma \ref{ylwithoutdelay}, we present the
following lemma.

\begin{yl}\label{yldelay}Consider the asynchronous case of system
\eqref{xtdelay}. For any $i\in\n, k\in\Z$, the number of elements
in set $\{{t}_j:t_j\in [t^i_{k}, t^i_{k+1})\}$ is not greater than
$\mml n (K(n-1)+1)$. We let $\mmm=\mml n (K(n-1)+1)$ and
$\mmu=(K+1)\mmm$;
\end{yl}

Proof: Obviously, there can not be infinite elements in $\{t_j:
t_j\in [t^i_{k}, t^i_{k+1})\}$ by Assumption (A1) and (A3). We
work out one upper bound of the number of elements in this set.

For any $i$, let $\{t_k^{(i)}\}=\{t_k^i-\tau_{ij}^k,k\in\Z, j\in
\N(t_k^i,i)\}\cup\{t_k^i\}$. Given $k\in\Z,i\in\n$, by Assumption
(A3), all possible elements in
$\{t_k^{(i)}\}\cap(t_k^i,t_{k+1}^i)$ are
$t_{k+1}^i-\tau_{ij}^{k+1}, j\in \N(t_{k+1}^i,i)$,
$t_{k+2}^i-\tau_{ij}^{k+2}, j\in \N(t_{k+2}^i,i)$, $\cdots$,
$t_{k+K}^i-\tau_{ij}^{k+K}, j\in \N(t_{k+K}^i,i)$. Because
$|\N(t_{k'}^i,i)|\leq n-1$ for any $k'$, where $|\N(t_{k'}^i,i)|$
is the number of elements in set $\N(t_{k'}^i,i)$. Therefore
\begin{equation}\label{eqyldelay}
    |\{t_k^{(i)}\}\cap[t_k^i,t_{k+1}^i)|\leq K(n-1)+1.
\end{equation}

Let $i\in\n$, $k\in\Z$ be given. By Lemma \ref{ylwithoutdelay},
$|\{t_{k'}^j:t_{k'}^j \in[t_k^i,t_{k+1}^i], k'\in \Z,j\in\n\}|\leq
\mml+1$. Let its elements be
$\bar{t}_1,\bar{t}_2,\cdots,\bar{t}_{m'+1}$ such that
$\bar{t}_{k'}<\bar{t}_{k'+1}$, $k'=1,2,\cdots, m'$, where
$m'\leq\mml$. For any $k'\in\{1,2,\cdots,m'\}$, $j\in\n$, there
exists $k''$ such that $t_{k''}^j\leq
\bar{t}_{k'}<\bar{t}_{k'+1}\leq t_{k''+1}^j$. From
\eqref{eqyldelay},
\[|[\bar{t}_{k'},\bar{t}_{k'+1})\cap\{t_k^{(j)}\}|\leq|[t^j_{k''},t^j_{k''+1})\cap\{t_k^{(j)}\}|\leq
K(n+1)+1.\] Since
$\{t_j:t_j\in[t_k^i,t_{k+1}^i)\}=\bigcup_{k'=1}^{m'}\left([\bar{t}_{k'},
\bar{t}_{k'+1})\cap\left(\cup_{j=1}^n\{t_k^{(j)}\}\right)\right)$,
  $|\{t_j: t_j\in[t_k^i,t_{k+1}^i)\}|\leq \mml n(K(n-1)+1)$.

For agent $i$ and any $k\in\Z$, $k \geq \mmu-1$, there exists $s$
 such that $t^i_s
\leq t_k< t_{k+1}\leq t^i_{s+1}$, and  agent $i$'s
 dynamics can be written
\[
    \dot{x}_i(t)=\left\{
                   \begin{array}{ll}
                     \sum_{j\in\N(t_s^i,i)}\sa_{ij}(t_s^i)(x_j(t_s^i-\tau^s_{ij})-x_i(t)), & \mbox{if\ }\N(t_s^i,i)\not=\phi; \\
                     0, & \mbox{otherwise}, \\
                   \end{array}
                 \right.
\]
where $t\in[t_k, t_{k+1})$. From Lemma \ref{yldelay} and
Assumption (A3), $t_{k-\mmu+1}\leq t_s^i-\tau_{ij}^s$. Let
$\tau_{ii}^s=0$. Solving the above equation gives
\[x_i(t_{k+1})=e^{-\tau_k}
x_i(t_{k})+(1-e^{-\tau_k})\sum_{j\in\N(t_s^i,i)\cup\{i\}}\sa_{ij}(t_k)x_j(t^i_s-\tau^s_{ij}).\]
 Let
$z(k)=[x(t_k)^T,x(t_{k-1})^T,\cdots,x(t_{k-\mmu+1})^T]^T$. For any
$k\geq \mmu-1$, there exists a matrix
$\pi(\tau_k,\sA_{1}(t_k),\cdots,\sA_{{\mmu}}(t_k))\in\Pi(\mmu,t_k)$
such that
\begin{equation}\label{xtdelay2}
    z(k+1)=\pi(\tau_k,\sA_{1}(t_k),\cdots,\sA_{{\mmu}}(t_k))z(k).
\end{equation}

With the same arguments as Proposition \ref{mt1}, we have
\begin{mt}\label{mt2}
System \eqref{xtdelay} solves a consensus problem if and only if
system \eqref{xtdelay2} solves a consensus problem.
\end{mt}

Now, we present the main result of this paper.
\begin{dl}\label{dlmain}
If there exists $T\geq 0$ such that for all $t^0\geq 0$, the union
of graph $\Ge$ across interval $[t^0, t^0+T]$ contains a spanning
tree, then the time-delayed system \eqref{xtdelay} solves a
consensus problem. Moreover, if $\G$ is time-invariant, then the
solvability of the consensus problem of system \eqref{xtdelay}
implies that $\G$ contains a spanning tree.
\end{dl}

\begin{rmk}
The second part of Theorem \ref{dlmain} is obvious.  In order to
prove Theorem \ref{dlmain}, it suffices to prove the possession of
consensus property of discrete-time system \eqref{xtdelay2} under
the hypothesis of Theorem \ref{dlmain}. However, the proof of
consensus property of system \eqref{xtdelay2} is not easy. The
following two properties of system \eqref{xtdelay2} make all
previous results, such as those in \cite{W. Ren and R. W. Beard,L.
Moreau 1}, inapplicable to our system, and therefore we have to
explore another efficient way.
\begin{enumerate}
    \item $\Pi(\mmu,t_k)$, the set including all possible state matrices $\pi(\tau_k,\cdot)$,
     is not a finite set, and furthermore is not  compact because that $\Gamma_{\ssA}$ is
     an
    infinite set and $\tau_k\in (0,\tuu]$, which is not compact.
    \item The diagonal entries of any $\pi(\tau_k, \cdot)$ are not all
    non-zeros.
\end{enumerate}
Since $\Pi(\mmu,t_k)$ is a infinite set, system \eqref{xtdelay2}
can be viewed as a discrete-time consensus model with the topology
that switches among an infinite set of directed graphs. Therefore,
previous results on discrete-time systems with finite available
topologies, such as those in \cite{A. Jadbabaie J. Lin and A. S.
Morse,R. Olfati-Saber and R. M. Murray 1,W. Ren and R. W. Beard},
are not applicable. The only general result touching on infinite
available topologies was given in \cite{L. Moreau 1}, while some
of the basic pre-requisitions such as the assumption of the
compactness of $e_k(\mathcal{A}(t))(x)$ and the assumption of the
strict convexity (namely, (3) of Assumption 1 in \cite{L. Moreau
1}) are not satisfied by our model (See Assumption 1 in \cite{L.
Moreau 1}).
\end{rmk}

The proof of Theorem \ref{dlmain}, presented in the next section,
is on the basis of the spacial structures of $\Pi(\mmu, t_k)$ and
properties of $\{t_k\}$.

\section{Technical Proof}
This section presents a complete proof of Theorem \ref{dlmain}. We
first give an equivalent formulation of the condition in Theorem
\ref{dlmain}.

\begin{yl}\label{ylconditionequi}
The existence of $T\geq 0$ such that for all $t^0\geq 0$, the
union of graph $\Ge$ across interval $[t^0, t^0+T]$ contains a
spanning tree, is equivalent to the condition that there exists
$\e\in\Z$ and $\tou>0$ with the following property:

For any $U_k=\{t_{k\e+1}, t_{k\e+2}, \cdots, t_{(k+1)\e}\}$, there
exists a subset of $U_k$, denoted by $V_k$, such that the union of
$\Ge$ on $V_k$ contains a spanning tree and for any $t_s\in V_k$,
$t_{s+1}-t_{s}\geq\tou$.
\end{yl}

Proof: The sufficiency is rather straightforward and only the
necessity is proved. By Lemma \ref{yldelay}, there exists $p\in\Z$
such that for any $k\in\Z$, $t_{k+p}-t_k\geq T$, such as
$p=\mmm(\lfloor\frac{T}{\tul}\rfloor+2)$ (For any given $i\in\n$,
there exists $k'$ such that $t_{k'}^i\leq t_k<t_{k'+1}^i$, and
$t^i_{k'+\frac{p}{\mmm}}-t_{k'+1}^i\geq\tul(\frac{p}{\mmm}-1)=\tul(\lfloor\frac{T}{\tul}\rfloor+1)>T$.
We claim that $t_{k+p}\geq t_{k'+\frac{p}{\mmu}}^i$. If not,
$|[t_{k'}^i,t^i_{k'+\frac{p}{\mmm}})\cap\{t_k\}|\geq|[t_k,t_{k+p}]\cap\{t_k\}|=p+1$.
But by Lemma \ref{yldelay},
$|[t_{k'}^i,t_{k'+\frac{p}{\mmm}}^i)\cap\{t_k\}|\leq
\frac{p}{\mmm}\mmm=p$, which is a contradiction. Therefore
$t_{k+p}-t_k\geq t_{k'+\frac{p}{\mmm}}^i-t_k>
t_{k'+\frac{p}{\mmm}}^i-t_{k'+1}^i>T$). Let $\e=p+ 2\mmm$. We
consider $U_k$. Obviously $t_{k\e+\mmm+p}-t_{k\e+\mmm}\geq T$.
Therefore the union graph $\Ge$ on
$\{t_{k\e+\mmm},\cdots,t_{k\e+\mmm+p}\}$ contains a spanning tree.
Let the edge set of the spanning tree  be $\mathcal{E}$. If
$(v_j,v_i)\in\mathcal{E}$, there exists $k'$ such that
$k\e+\mmm\leq k'\leq k\e+\mmm+p$ and
$(v_j,v_i)\in\mathcal{G}^0(t_{k'})$. For $t_{k'}$, there exists
$k''$ such that $t^i_{k''}\leq t_{k'}< t^i_{k''+1}$, and thus
$(v_j,v_i)$ is an edge of the graph $\Ge$, $t\in[t_{k''}^i,
t^i_{k''+1})$.

We claim that $t_{k''}^i\geq t_{k\e+1}$ and $t_{k''+1}^i\leq
t_{(k+1)\e}$. In fact, if $t_{k''}^i<t_{k\e+1}$, then
$|[t_{k''}^i,t_{k''+1}^i)\cap\{t_k\}|
>|[t_{k\e+1},t_{k'}]\cap\{t_k\}|\geq|[t_{k\e+1},t_{k\e+\mmm}]\cap\{t_k\}|=\mmm$,
which contradicts Lemma \ref{yldelay}. And if
$t_{k''+1}^i>t_{(k+1)\e}$, then
$|[t_{k''}^i,t_{k''+1}^i)\cap\{t_k\}|
\geq|[t_{k'},t_{(k+1)\e}]\cap\{t_k\}|\geq|[t_{k\e+\mmm+p},t_{(k+1)\e}]\cap\{t_k\}|=\mmm+1$,
which also  contradicts Lemma \ref{yldelay}.

 Since $\{t_k\}\cap[t^i_{k''},t_{k''+1}^i)$ has at
most $\mmm$ elements, there exists $t_s\in
\{t_k\}\cap[t^i_{k''},t_{k''+1}^i)\subset U_k$ such that
$t_{s+1}-t_s\geq\frac{t_{k''+1}^i-t_{k''}^i}{\mmm}\geq\frac{\tul}{\mmm}$.
Let $\tou=\frac{\tul}{\mmm}$. If $(v_j,v_i)$ takes very possible
edge in $\mathcal{E}$, we obtain all possible $t_s$s. Let the set
of them be $V_k$. Then $V_k$ has the aforementioned property and
the necessity is proved.

Let $A, B$ be  $r\times r$ stochastic matrices and let
$\delta(A)=\max_j\max_{i_1,i_2}|a_{i_1j}-a_{i_2j}|$. Thus
$\delta(A)$ measures how different the rows of $A$ are. If the
rows of $A$ are identical, $\delta(A)=0$ and conversely. We say
that $A, B$ are of the same type, $A\sim B$, if they have zero
elements and positive elements in the same place. Let $\num(r)$ be
the number of different types of all $r\times r$ SIA matrices.
Define $\lambda(A)=1-\min_{i_1,i_2}\sum_j\min(a_{i_1,j},
a_{i_2,j})$. If $\lambda(A)<1$ we call $A$ a {\it scrambling
matrix}.

\begin{yl}[\cite{J. Wolfowitz}, Lemma 2]\label{yldelta}
For any stochastic matrices $A_1, A_2, \cdots, A_k$, $k>0$,
\[\delta(A_1 A_2\cdots A_k)\leq \prod_{i=1}^k\lambda(A_i).\]
\end{yl}

The next lemma generalizes the result of Lemma 4 in \cite{J.
Wolfowitz}.
\begin{yl}\label{ylscrambling}
Let $A_1, A_2,\cdots, A_k$ (repetitions permitted) be $r \times r$
SIA matrices with the property that for any $1\leq k_1< k_2\leq
k$, $\prod_{i=k_1}^{k_2} A_i$ is SIA. If $k>\num(r)$, then
$\prod_{i=1}^k A_i$ is a scrambling matrix.
\end{yl}
Proof: Since $k>\num(r)$, there exist $k_1<k_2$, such that
$\prod_{i=1}^{k_1}A_i\sim \prod_{i=1}^{k_2}A_i$. It follows from
Lemma \ref{ylWolfowitz3} and $\prod_{i=k_1+1}^{k_2}$ being SIA
that $\prod_{i=1}^{k_1}A_i$ is a scrambling matrix. Thus
$\prod_{i=1}^kA_i$ is also a scrambling matrix by Lemma
\ref{ylWolfowitz1}.

To investigate the properties of matrices in $\Pi(\mmu, t_k)$, we
introduce some notations. Let $\F(A)=\sum_{i=1}^{\mmu}A_{1i}$,
where $A=[A_{ij}]$ is an $\mmu\times\mmu$ block matrix and
$A_{ij}\in\R^{n\times n}$. Let $\Gamma_s$ denote the set of square
matrices such that $A\in\Gamma_s$ if and only if $\mathcal{G}(A)$
contains a spanning tree with the property that the root vertex of
the spanning tree has a self-loop in $\mathcal{G}(A)$.

\begin{yl}\label{ylmain}
Let $A$ be a stochastic matrix. If $A\in\Gamma_s$, then $A$ is
SIA.
\end{yl}
Proof: Since $A$ is stochastic, $A\yi=\yi$ and $\rho(A)=1$. We
assume that there exists a spanning tree with vertex $v_{k_1}$ as
its root and $(v_{k_1}, v_{k_1})\in \mathcal{E(G}(A))$. Suppose
that subgraph $\mathcal{G}_s$ induced by $v_{k_1},
v_{k_2},\cdots,v_{k_s}$ ($1\leq s \leq n$) is the maximal induced
subgraph that is strongly connected. Let the vertices in
$\mathcal{V}(\mathcal{G}(A)) \backslash \{v_{k_1},
v_{k_2},\cdots,v_{k_s}\}$ be $v_{k_{s+1}} ,\cdots,v_{k_n}$. Then
there exists a permutation matrix $T$ such that \[\left[
         \begin{array}{c}
           k_1 \\
           k_2 \\
           \vdots \\
           k_n \\
         \end{array}
       \right]=T\left[
         \begin{array}{c}
           1 \\
           2 \\
           \vdots \\
           n \\
         \end{array}
       \right].
\]
Therefore,
\[
TAT^{-1}=\left[
           \begin{array}{cc}
             A_{11} & A_{12} \\
             A_{21} & A_{22} \\
           \end{array}
         \right],
\]
where $A_{11}=\left[
                \begin{array}{ccc}
                  a_{k_1,k_1} & \cdots & a_{k_1,k_s} \\
                  \vdots& \ddots & \vdots \\
                  a_{k_s,k_1} &\cdots & a_{k_s,k_s} \\
                \end{array}
              \right]
$, $A_{12}=\left[
                \begin{array}{ccc}
                  a_{k_{1},k_{s+1}} & \cdots & a_{k_1,k_n} \\
                  \vdots& \ddots & \vdots \\
                  a_{k_s,l_{s+1}} &\cdots & a_{k_s,k_n} \\
                \end{array}
              \right]
$, \break $A_{21}=\left[
                \begin{array}{ccc}
                  a_{k_{s+1},k_1} & \cdots & a_{k_{s+1},k_s} \\
                  \vdots& \ddots & \vdots \\
                  a_{k_n,k_1} &\cdots & a_{k_n,k_s} \\
                \end{array}
              \right]
$ and $A_{22}=\left[
                \begin{array}{ccc}
                  a_{k_{s+1},k_{s+1}} & \cdots & a_{k_{s+1},k_n} \\
                  \vdots& \ddots & \vdots \\
                  a_{k_n,k_{s+1}} &\cdots & a_{k_n,k_n} \\
                \end{array}
              \right].
$

By the assumption that $\mathcal{G}_s$ is maximal, $A_{12}=0$.
Since $\mathcal{G}_s$ is strongly connected, $A_{11}$ is
irreducible. And from
$(v_{k_1},v_{k_1})\in\mathcal{E}(\mathcal{G}(A))$,
$a_{k_1,k_1}>0$. By Lemma \ref{yl848}, $A_{11}$ is primitive (Its
definition is provided in the Appendix) and thus $1$ is the only
eigenvalue of $A_{11}$ with maximum modulus. Since $1$ is an
eigenvalue of $A_{11}$, by Lemma \ref{RWyl}, $1$ is not an
eigenvalue of $A_{22}$. On the other hand, let $\rho(A_{22})$
denote the spectral radius of $A_{22}$. By Ger\v sgorin disc
theorem, $\rho(A_{22})\leq 1$ and by Lemma \ref{yl831}, $\rho(A)$
is an eigenvalue of $A$. It follows that $\rho(A_{22})<1$.
Consequently, $1$ is the only eigenvalue of $A$ with maximum
modulus. By Lemma \ref{RWyl} and that $\rho(A)=1$, it is easy to
check that $A$ satisfies the conditions of Lemma \ref{yl827}. Let
$f^TA=f^T$ such that $f^T\yi=1$. Then $\lim_{k\to\infty}A^k=\yi
f^T$.

\begin{figure}[htpb]
\centering
        \includegraphics[scale=0.6]{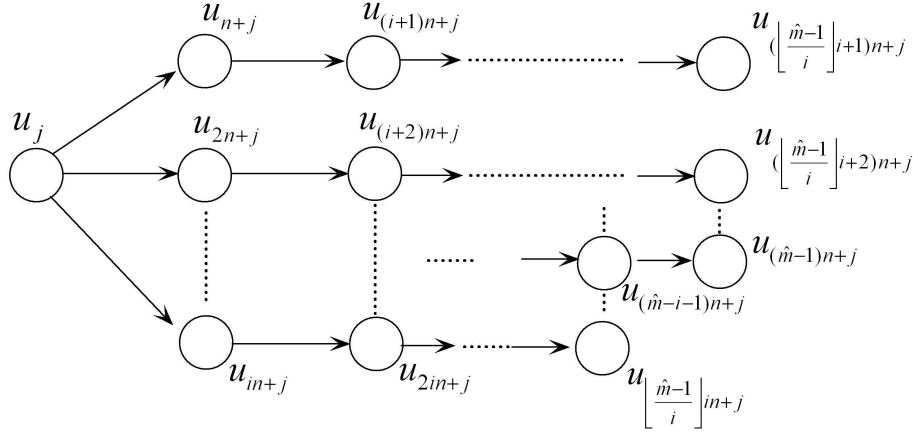}
      \caption{Paths starting from $u_j$ in $\mathcal{G}(M_0^i)$. If $\lfloor\frac{\hat{m}-1}{i}\rfloor i=\hat{m}-1$, $u_{\lfloor\frac{\hat{m}-1}{i}\rfloor in+j}$
       and $u_{(\hat{m}-1)n+j}$ are the same vertex.}
      \label{fig1}
\end{figure}

\begin{yl}{\label{yltu}}
Let $A_1,\cdots,A_{\mmu}$ be $n \times n$ nonnegative matrices and
let
\begin{equation}
\left.
  \begin{array}{c}
   M_0=\left[
        \begin{array}{ccccc}
          I &   &   &   &   \\
          I &  & &\mbox{\Huge 0 } &  \\
           & I & &&  \\
           &  & \ddots &  &  \\
          & \mbox{\Huge 0 }&  & I & 0 \\
        \end{array}
      \right]_{\mmu n\times \mmu n}
, M_1=\left[
        \begin{array}{ccccc}
          I+A_1 & A_2 & \cdots & A_{\mmu-1} & A_{\mmu} \\
          I &  & & &  \\
           & I & &\mbox{\Huge 0 }&  \\
           &  & \ddots &  &  \\
          & \mbox{\Huge 0 }&  & I &  \\
        \end{array}
      \right],\\
M_2=\left[
        \begin{array}{ccccc}
          I +A_1& \cdots & A_{m-2} & A_{\mmu-1} & A_{\mmu} \\
          I &  &  &  & \\
          I &  &  &\mbox{\Huge 0 }  &  \\
           & \ddots & & & \\
           \multicolumn{2}{c}{\mbox{\Huge 0 }} & I &  &  \\
        \end{array}
      \right],
    \cdots,M_{\mmu-1}=\left[
        \begin{array}{ccccc}
          I+A_1 & A_2 &A_3& \cdots  & A_{\mmu} \\
          I &  &  &  &\\
          I &  & \multicolumn{2}{c}{\raisebox{-2ex}[0pt]{\mbox{\Huge0}}}  &\\
          \vdots & & & &  \\
          I &  &  & & \\
        \end{array}
      \right].\\
  \end{array}
\right.
\end{equation}
For any $i\in\{1,2,\cdots,\mmu-1\}$, if $\mathcal{G}(\F(M_i))$
contains a spanning tree, then $M_i\in\Gamma_s$.
\end{yl}
Proof: Let $N=M_1-M_0$. Then we have $M_i=M_0^i+N$. Let
$\mathcal{G}(M_i)$, $\mathcal{G}(M_0^i)$, and $\mathcal{G}(N)$ be
with the same vertex set $\{u_1,u_2,\cdots,u_{\mmu n}\}$ and let
$\mathcal{E}(\mathcal{G}(\F(M_i)))=\{v_1,v_2,\cdots,v_n\}$.
Apparently
$\mathcal{E}(\mathcal{G}(M_i))=\mathcal{E}(\mathcal{G}(M_0^i))\cup\mathcal{E}(\mathcal{G}(N))$.
We first investigate the edge sets
$\mathcal{E}(\mathcal{G}(M_0^i))$ and
$\mathcal{E}(\mathcal{G}(N))$. For any $j\in\n$, $0\leq i \leq
\mmu-1$,
\begin{align*}
    &\{(u_j,u_{n+j}),(u_j, u_{2n+j}),\cdots,  (u_j,u_{in+j}),\\
    &(u_{n+j}, u_{(i+1)n+j}), (u_{2n+j},u_{(i+2)n+j}),\cdots,
    (u_{(\mmu-i-1)n+j},u_{(\mmu-1)n+j})\}\subset\mathcal{E}(\mathcal{G}(M_0^i)).
\end{align*}
Therefore for any $j\in\n$, there exist paths  from $u_j$ to
$u_{n+j},\cdots, u_{(\mmu-1)n+j}$ in
$\mathcal{E}(\mathcal{G}(M_0^i))$ (See Fig. \ref{fig1}). If there
exists an edge $(v_j,v_k)\in\mathcal{E}(\mathcal{G}(\F(M_i)))$,
then there exists $0\leq s\leq \mmu-1$ such that
$(u_{j+sn},u_k)\in\mathcal{E}(\mathcal{G}(N))$. Therefore there
exists a path from $u_j$ to $u_k$ in $\mathcal{G}(M_i)$. It
follows that if $\mathcal{G}(\F(M_i))$ contains a spanning tree
with root vertex $v_j$, $\mathcal{G}(M_i)$ also contains a
spanning tree with root vertex $u_j$. Since the entry in the $j$th
row and the $j$th column of $M_i$ is not less than $1$, $u_j$ has
a self-loop in $\mathcal{G}(M_i)$. Consequently
$\mathcal{G}(M_i)\in\Gamma_s$.
Let $\Pi(\mmu)$  be the set of matrices \begin{equation*} \left[
                            \begin{array}{ccccc}
                              e^{-h}I+(1-e^{-h})A_1 & (1-e^{-h})A_2 & \cdots & (1-e^{-h})A_{m-1} & (1-e^{-h})A_{\mmu} \\
                              I & 0 & \cdots & 0 & 0 \\
                               0 & I& \cdots & 0 & 0 \\
                              \vdots & \vdots & \ddots & \vdots & \vdots \\
                              0 & 0 & \cdots & I &0 \\
                            \end{array}
                          \right]_{\mmu n\times \mmu n},
\end{equation*}
where $0\leq h\leq \tuu$, and there exists some
$\sA_0\in\Gamma_{\ssA}$, such that $
A_1,\cdots,A_{\mmu}\in\Lambda(\sA_0)$, and $
A_1+\cdots+A_{\mmu}=\sA_0$. Since $\Gamma_{\ssA}$ is compact and
given any $\sA_0$, all possible choices of $A_1,\cdots,A_{\mmu}$
are finite, $\Pi(\mmu)$ is a compact set, and for any $k\in\Z$,
$\Pi(\mmu,t_k)\subset\Pi(\mmu)$.

Let $\Pi_0=\{\Pi_{i=1}^{\e}\pi(h_i,
A_{i1},A_{i2},\cdots,A_{i\mmu}): \pi(h_i,\cdot)\in \Pi(\mmu),$ and
there exists a subset of $\{1,2,\cdots,{\e}\}$, denoted by
$\mathcal{H}$, such that for any $s\in\mathcal{H}$, $h_s\geq\tou$,
and $\mathcal{G}(\sum_{i\in\mathcal{H}}\sum_{j=1}^{\mmu}A_{ij})$
contains a spanning tree$\}$, where $\e$ and $\tou$ are defined in
Lemma \ref{ylconditionequi}.

\begin{yl}
$\Pi_0$ is a compact set and for any $\pi\in\Pi_0$, $\pi$ is SIA.
And for any $k>\num(n\mmu)$,   if $\pi_1$, $\cdots$,
$\pi_{k}\in\Pi_0$, then $\prod_{i=1}^k\pi_i$ is a scrambling
matrix and there exists a $\hat{\lambda}(k)\in[0,1)$, such that
$\lambda(\prod_{i=1}^k\pi_{i})\leq \hat{\lambda}(k)$.
\end{yl}
Proof: The compactness of $\Pi_0$ follows from the following facts
\begin{enumerate}
    \item $ \Pi(\mmu)$ is a compact set;
    \item All possible choices of $\mathcal{H}$ are finite;
    \item $h_s\in[\tou, \tuu]$, which is a compact set;
    \item All possible choices of the spanning tree are finite;
    \item Given the spanning tree and $\mathcal{H}$, let
    $\Pi_1=\{\prod_{i=1}^{\e}\pi(h_i, A_{i1},\cdots, A_{i\mmu}):\pi(h_i,\cdot)\in\Pi(\mmu)$,
    and for any $s\in\mathcal{H}$, $h_s\in[\tou,\tuu]$, and
    $\mathcal{G}(\sum_{i\in\mathcal{H}}\sum_{j=1}^{\mmu}A_{ij})$
    contains the specified spanning tree$\}$ is compact.
\end{enumerate}

We only prove the fact 5)

 Let $|\mathcal{H}|=q$, $q\leq\e$. Since the
 product of $\e$ matrices is continuous, it suffices to
 prove that $\Pi_2=\{[\pi(h_1, A_{11},\cdots, A_{1\mmu}), \cdots,\pi(h_{\e}, A_{\e1},\cdots, A_{\e\mmu})]:
 \pi(h_i,\cdot)\in\Pi(\mmu)$,
    and for any $s\in\mathcal{H}$, $h_s\in[\tou,\tuu]$, and
    $\mathcal{G}(\sum_{i\in\mathcal{H}}\sum_{j=1}^{\mmu}A_{ij})$
    contains the specified spanning tree
    $\}$ is compact. Since $\pi(h_i,\cdot)\in \Pi(\mmu)$, which
    is compact,  it suffices to prove that
    $\Pi_3=\{[\pi(h_{1}, A_{1 1},\cdots, A_{1\mmu}), \cdots,\pi(h_{q}, A_{q},$ $\cdots, A_{q\mmu})]:
 \pi(h_{i},\cdot)\in\Pi(\mmu)$,
     $h_i\in[\tou,\tuu]$, $i=1,2,\cdots,q$, and
    $\mathcal{G}(\sum_{i=1}^q\sum_{j=1}^{\mmu}A_{ij})$
    contains the specified spanning tree
    $\}$ is compact.

Let $B^{(p)}=[b^{(p)}_{ij}]\in\Pi_3, p=1,2,\cdots$ be a sequence
of matrices , and $\lim_{p\to\infty}B^{(p)}=B=[b_{ij}]$. Since
$\Pi_3$ is a bounded set, it suffices to prove that $B\in\Pi_3$.

It is clear that $\lim_{p\to\infty}b^{(p)}_{ij}=b_{ij}$. For any
$i,j\in\n$, if $b_{ij}^{(p)}\not=0$, then by the definition of
$\Gamma_{\ssA}$,
$b_{ij}^{(p)}>(1-e^{-\tou})\frac{\mathbf{\check{a}}}{(n-1)\mathbf{\hat{a}}}$.
Therefore if $b_{ij}>0$, then
$\lim_{p\to\infty}b_{ij}^{(p)}=b_{ij}\geq
(1-e^{-\tou})\frac{\mathbf{\check{a}}}{(n-1)\mathbf{\hat{a}}}>0$,
and thus there exists $P_{ij}\in\Z$ such that for any $p>P_{ij}$,
$b^{(p)}_{ij}> 0$. If $b_{ij}=0$, from
$\lim_{p\to\infty}b_{ij}^{(p)}=b_{ij}=0$, there exists
$P_{ij}\in\Z$ such that for any $p>P_{ij}$, $b^{(p)}_{ij}=0$. Let
$P=\max_{ij}P_{ij}$, and then for any $p>P$, $B^{(p)}\sim B$.

Let $B=[B_1, B_2,\cdots,B_q]$, $B_i\in\R^{\mmu n\times\mmu n}$,
$1\leq i\leq q$. From the compactness of $\Pi(\mmu)$,
$B_i\in\Pi(\mmu)$. If $B_i=\pi(h_{b_i},B_{i1},\cdots,B_{i\mmu})$,
$1\leq i\leq q$, then $h_{b_i}\in[\tou,\tuu]$. And from
$B^{(p)}\sim B$,
$\mathcal{G}(\sum_{i=1}^q\sum_{j=1}^{\mmu}B_{ij})$ contains the
spanning tree. To conclude, $B\in\Pi_3$.

Next, we prove that for any $\pi\in\Pi_0$, $\pi$ is SIA. Let
$\pi=\sum_{i=1}^{\e}\pi(h_i, A_{i1}, \cdots,A_{i\mmu})$ and
$\mathcal{H}$ be the associated subset of $\{1,2,\cdots,\e\}$
defined in the definition of $\Pi_0$. Let $M_0$ be the same as in
Lemma \ref{yltu} and let
\[D_i=\left[
       \begin{array}{cccc}
         A_{i1} & A_{i2} & \cdots & A_{i\mmu} \\
         0 & 0 & \cdots & 0 \\
         \vdots & \vdots & \ddots & \vdots \\
         0 & 0 & \cdots & 0\\
       \end{array}
     \right].
\]
\begin{align*}
    \prod_{i=1}^{\e}\pi(h_i,\cdot)\geq & \prod_{i=1}^{\e}
(e^{-h_i}M_0+(1-e^{-h_i})D_i)\\
\geq
&e^{-\e\tuu}M_0^{\e}+e^{-(\e-1)\tuu}\sum_{i=1}^{\e}(1-e^{-h_i})(M_0)^{i-1}D_i(M_0)^{\e-i}\\
\geq &
e^{-\e\tuu}M_0^{\e}+e^{-(\e-1)\tuu}\sum_{i\in\mathcal{H}}(1-e^{-h_i})D_iM_0^{\e-i}\\
\geq&e^{-\e\tuu}M_0^{\e}+e^{-(\e-1)\tuu}(1-e^{-\tou})\sum_{i\in\mathcal{H}}D_iM_0^{\e-i}\\
\geq &
\min\{e^{-\e\tuu},e^{-(\e-1)\tuu}(1-e^{-\tou})\}(M_0^{\e}+\sum_{i\in\mathcal{H}}D_iM_0^{\e-i}).
\end{align*}
The second inequality follows from $0< h_i\leq \tuu$, the third
follows from $M_0^{i-1}D_i\geq D_i$, and the fourth follows from
$h_i\geq \tou, i\in\mathcal{H}$.

From $\F(D_iM_0^{\e-i})=\F(D_i)$, we have
$\F(\sum_{i\in\mathcal{H}}D_iM_0^{\e-i})=\sum_{i\in\mathcal{H}}\F(D_iM_0^{\e-i})=\sum_{i\in\mathcal{H}}\F(D_i)=\sum_{i\in\mathcal{H}}\sum_{j=1}^{\mmu}A_{ij}$.
Since $\mathcal{G}(\sum_{i\in\mathcal{H}}\sum_{j=1}^{\mmu}A_{ij})$
contains a spanning tree,
$\mathcal{G}(\F(\sum_{i\in\mathcal{H}}D_iM_0^{\e-i}))$ also
contains a spanning tree. Let $N\in\R^{\mmu n\times\mmu n}$ be
with the same first $n$ rows as
$\sum_{i\in\mathcal{H}}D_iM_0^{\e-i}$ and all other rows are
zeros. Then $ \prod_{i=1}^{\e}\pi(h_i,\cdot)\geq
\min\{e^{-\e\tuu},e^{-(\e-1)\tuu}(1-e^{-\tou})\}(M_0^{\e}+N)$. By
Lemma \ref{yltu}, $M_0^{\e}+N\in\Gamma_s$, and thus
$\pi\in\Gamma_s$. $\pi$ is also stochastic, and therefore, by
Lemma \ref{ylmain}, $\pi$ is SIA.

With the same arguments, we can conclude that for any $1\leq
k_1<k_2\leq k$, $\prod_{i=k_1}^{k_2}\pi_i$ is SIA (We only need to
replace $\e$ by $(k_2-k_1+1)\e$ in the above arguments and
$\mathcal{H}$ be the index set associated to any $\pi_i$ as
defined in the definition of $\Pi_0$). By Lemma
\ref{ylscrambling}, for any $k>\num(n\mmu)$,
$\prod_{i=1}^k\pi_{i}$ is a scrambling matrix. Let
\[\hat{\lambda}(k)=\max_{\substack{\bar{\pi}_{i}\in\Pi_0\\1\leq i\leq k}}\lambda\left({\prod\bar{\pi}_{i}}\right).\]
Since $\Pi_{0}$ is a compact set and $\lambda(\cdot)$ is
continuous, $\hat{\lambda}(k)$ exists and $\hat{\lambda}(k)<1$.
Obviously $\lambda(\prod_{i=1}^k\pi_i)\leq\hat\lambda(k)$.

{\it Proof of Theorem \ref{dlmain}:}

For any $k\in\Z$, let
$\pi_k=\prod_{s=k\e+\mmu-1}^{(k+1)\e+\mmu-2}\pi(\tau_{s},\sA_1(t_s),
\sA_2(t_s),\cdots, \sA_{\mmu}(t_s))$. By Lemma
\ref{ylconditionequi} and Proposition \ref{mtgraphG},
$\pi_k\in\Pi_0$.

Let $p=\num(\mmu n)+1$. For any $q>p\e$, there exists $s\in\Z$
such that $q=sp\e+q'$, where $0\leq q'<p\e$. By Lemma
\ref{yldelta},
\[\delta(\prod_{k=\mmu-1}^{q+\mmu-2}\pi(\tau_k,\sA_{1}(t_k),\cdots,\sA_{\mmu}(t_k)))\leq \prod_{i=0}^{s-1}\lambda(\prod_{j=ip}^{(i+1)p-1}\pi_j)\leq (\hat{\lambda}(p))^s.\]
Therefore
\[ \lim_{q\to\infty}\delta(\prod_{k=\mmu-1}^{q+\mmu-2}\pi(\tau_k,\sA_{1}(t_k),\cdots,\sA_{\mmu}(t_k)))=0,\]
which implies that  there exists $a\in\R$ such that
$\lim_{q\to\infty}\prod_{k=\mmu-1}^q\pi(\tau_k,\sA_{1}(t_k),\cdots,\sA_{\mmu}(t_k))=\yi
a$. And thus system \eqref{xtdelay2} solves a consensus
problem.\hspace*{\fill}~
\par\endtrivlist\unskip

\section{Simulations}
In this section, we take some examples to illustrate the
effectiveness of our results.
\begin{lz}[Fixed topology without time-delays]\label{alz1}
Suppose that the system consists of $4$ agents. Let $\mathcal{A}=$
{\scriptsize $\left[\begin{array}{cccc}
                                                             0 & 1 & 1 & 0 \\
                                                             1 & 0 & 0 & 0 \\
                                                             0 & 1 & 0 & 0 \\
                                                             0 & 0 & 1 &0 \\
                                                           \end{array}
                                                         \right]$},
$\tul=0.2$ and $\tuu=0.9$. And suppose that there do not exist
communication time-delays and  each agent can get all its
neighbors' states at its update times.  For any $i\in\n$,
$t_{k+1}^i-t_k^i$ is evenly distributed between $0.2$ and $0.9$.
Since the communication topology has a spanning tree, the
consensus is reachable under asynchronous consensus control
strategy \eqref{xt1}. Let initial state $x(0)=[5, 6, 7, 8]^T$. In
the simulation experiment, the update times $t_k^i$, $i=1,2,3,4$,
are randomly generated and independent of each other. The state
trajectories of agents are shown in Fig. \ref{fig5}.

As a marked difference from the synchronous case, the final states
of agents are dependent on the update times. We repeat the
simulation experiment $100$ times independently, and the final
states of them are shown in Fig. \ref{fig6}.
\end{lz}
\begin{figure}[htpb]
\centering
        \includegraphics[scale=0.6]{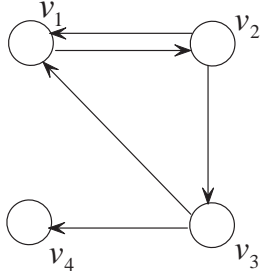}
      \caption{$\mathcal{G(A)}$ in Example \ref{alz1}}
      \label{fig4}
\end{figure}
\begin{figure}[htpb]
\centering
     \includegraphics[scale=0.6]{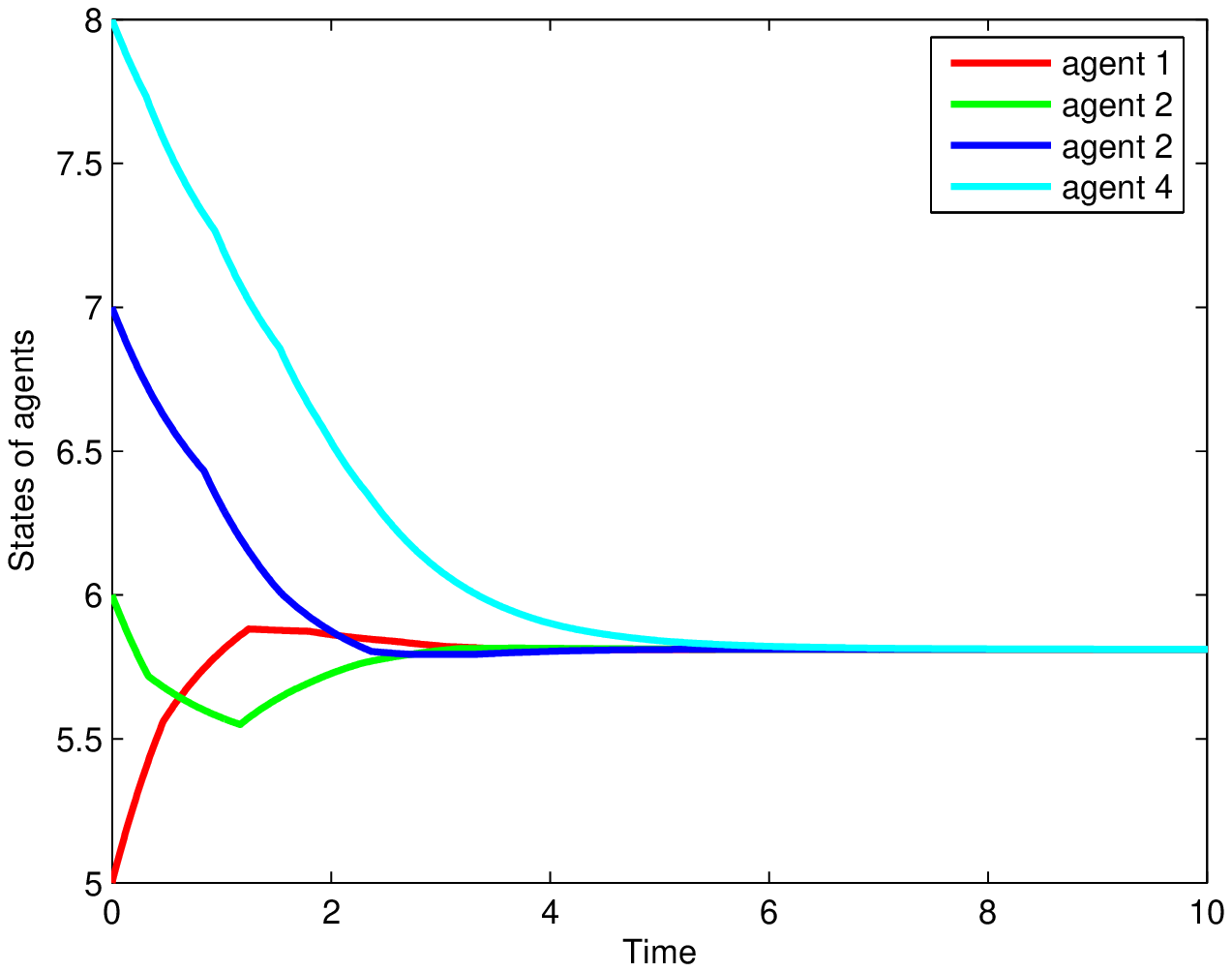}
      \caption{State trajectories of agents in Example \ref{alz1}}
      \label{fig5}
\end{figure}
\begin{figure}[htpb]
\centering
      \includegraphics[scale=0.6]{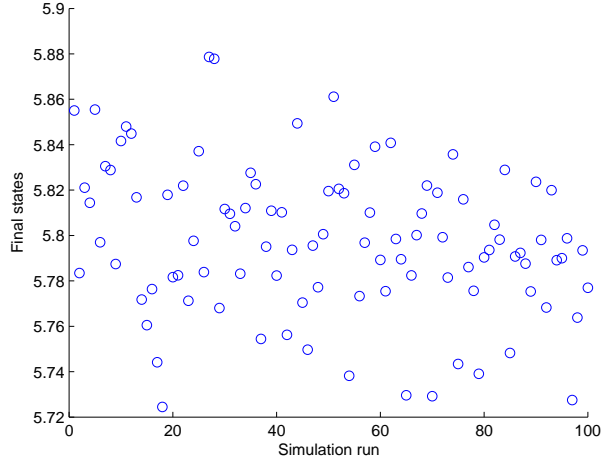}
      \caption{Final states of different experiments in Example \ref{alz1}}
      \label{fig6}
\end{figure}

\begin{figure}[htpb]
\centering        \includegraphics[scale=0.6]{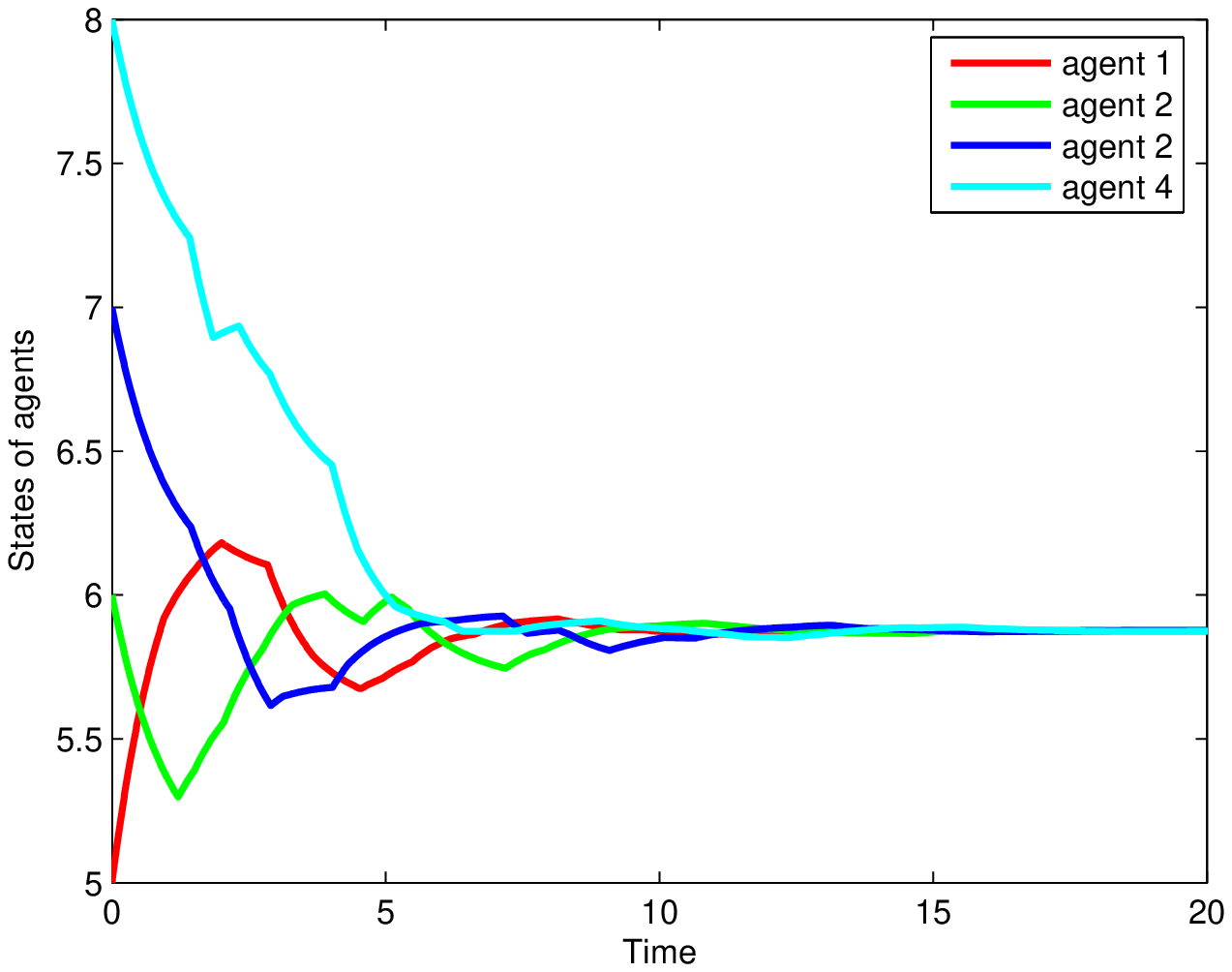}
      \caption{State trajectories of agents with $\tau_d=2$}
      \label{fig7}
\end{figure}
\begin{figure}[htpb]
\centering        \includegraphics[scale=0.6]{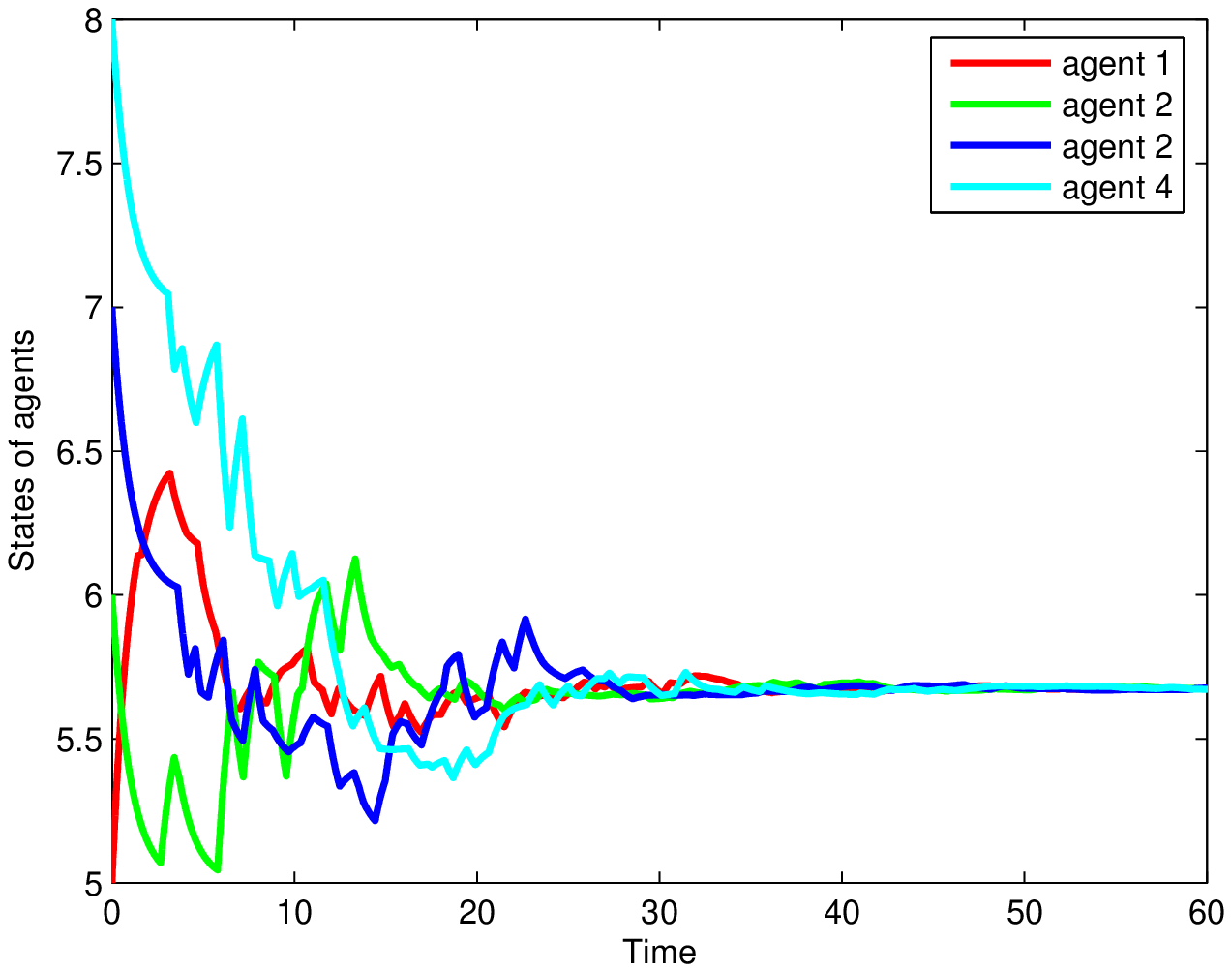}
      \caption{State trajectories of agents with $\tau_d=10$}
      \label{fig8}
\end{figure}\begin{figure}[htpb]
\centering
        \includegraphics[scale=0.6]{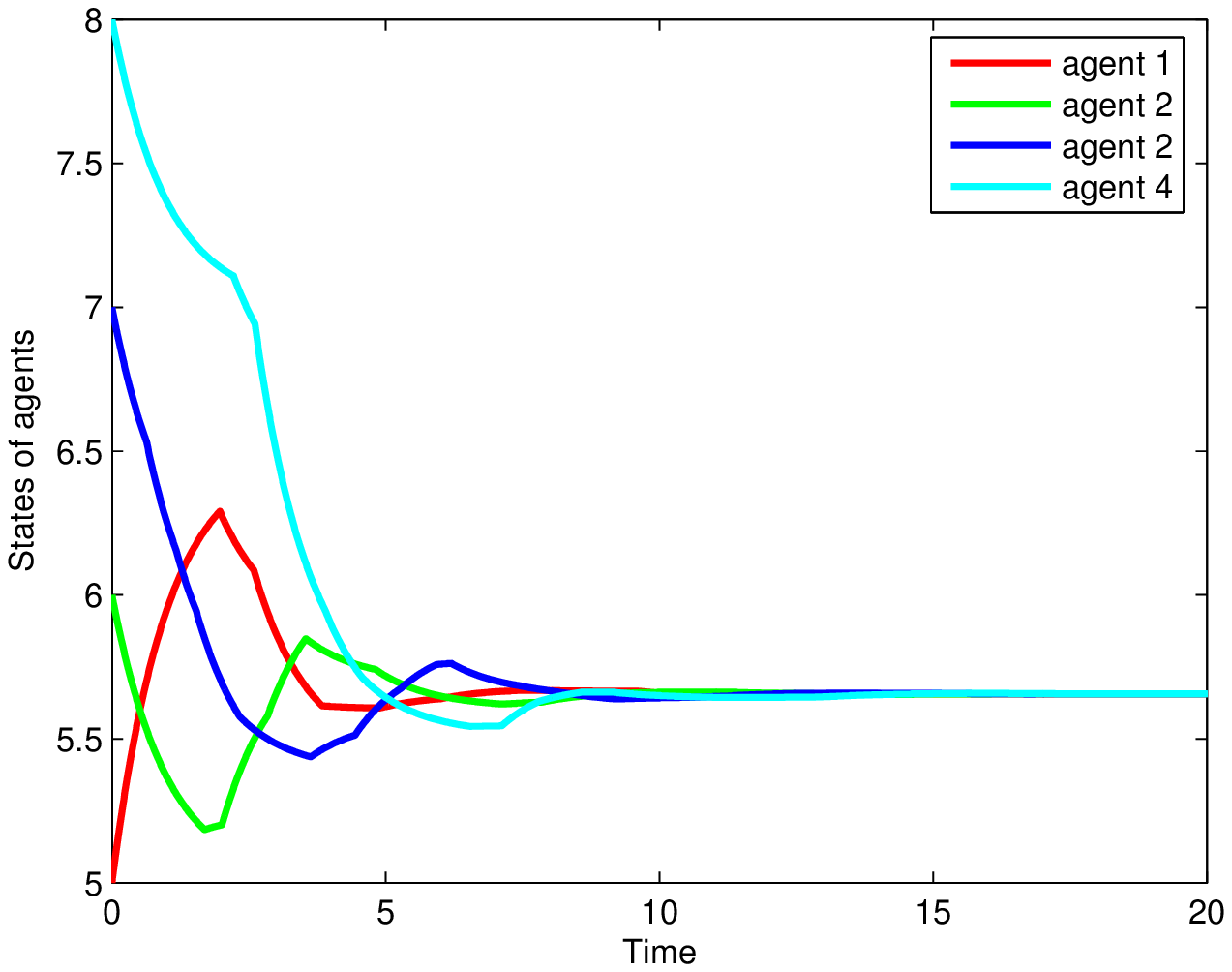}
      \caption{State trajectories of agents by the-most-recent-data strategy with $\tau_d=2$}
      \label{fig9}
\end{figure}\begin{figure}[htpb]
\centering
        \includegraphics[scale=0.6]{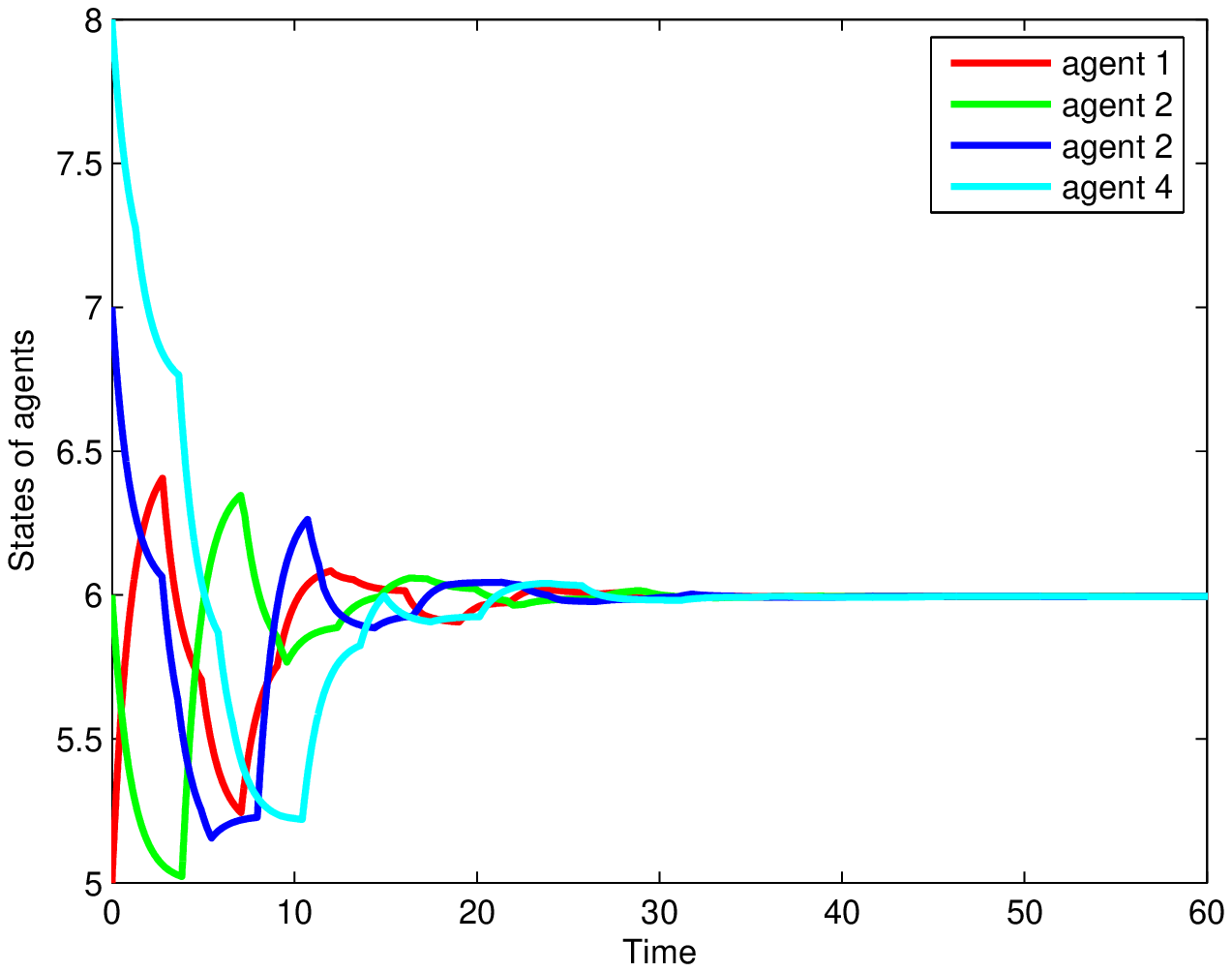}
      \caption{State trajectories of agents by the-most-recent-data strategy with $\tau_d=10$}
      \label{fig10}
\end{figure}
\begin{figure}[htpb]
\centering
        \includegraphics[scale=0.6]{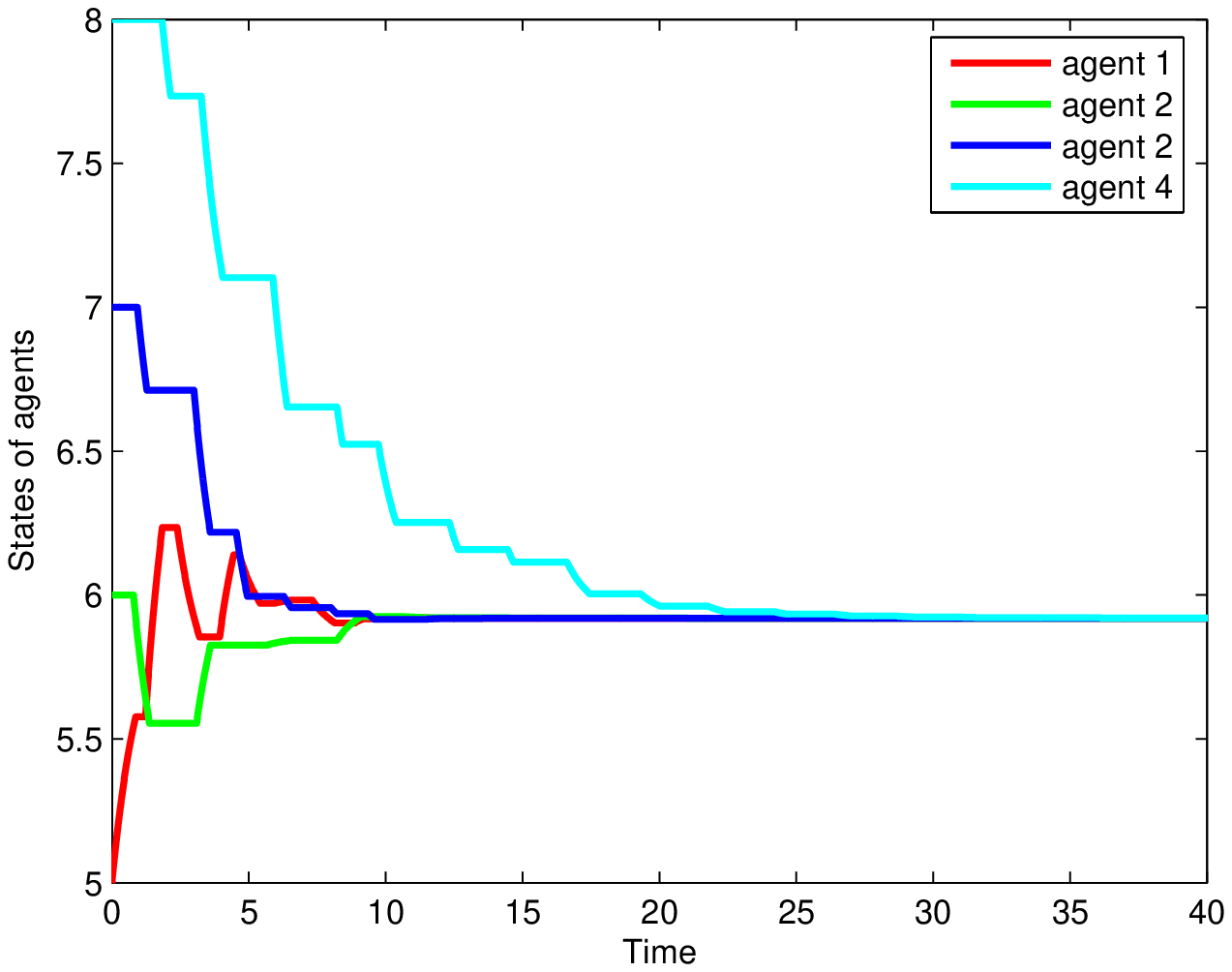}
      \caption{State trajectories of agents in Example \ref{lzswitching}}
      \label{switching}
\end{figure}

\begin{lz}[Fixed topology with time-delays]\label{lzdelay}
We still consider the system in Example \ref{alz1} and suppose
that there exist communication time-delays bounded by $\tau_d$ and
each agent can get all its neighbors' states at its update times.
We let initial state $x(t)=[5,6,7,8]$, $t\in[-\tau_d, 0]$. Fig.
\ref{fig7} and Fig. \ref{fig8} show the state trajectories of
agents with maximum communication time-delay $\tau_d=2$ ($K=10$)
and $\tau_d=10$ ($K=50$) separately, where the time-delays are
randomly generated. We can see that the system with $\tau_d=2$
converges faster than the system with $\tau_d=10$. If we adopt
the-most-recent-data strategy, we can get better convergence rate.
Fig. \ref{fig9} and Fig. \ref{fig10} show the state trajectories
of agents under the-most-recent-data strategy with randomly
generated time-delays bounded by $\tau_d=2$ and $\tau_d=10$
separately.
\end{lz}

\begin{lz}[Switching topology with time-delays]\label{lzswitching}
We  still consider the system consisting of $4$ agents. Each
agent's update intervals are evenly distributed between $0.2$ and
$0.9$, and are randomly generated. Suppose that the maximum
time-delay $\tau_d=2$, and initial state $x(t)=[5,6,7,8]$,
$t\in[-\tau_d,0]$.

We assume that the weight of each edge of the communication
topology is $1$ and
\begin{enumerate}
    \item agent $1$ can get the state of agent $2$ at update times
    $t^1_{4k}$, $k\in\Z$, and can get the state of agent $3$ at update
    times
    $t^1_{4k+2}$, $k\in\Z$;
    \item agent $2$ can get the state of agent $1$ at update times
    $t^2_{4k+1}$, $k\in\Z$;
    \item agent $3$ can get the state of agent $2$ at update times
    $t^3_{4k+2}$, $k\in\Z$;
    \item agent $4$ can get the state of agent $3$ at update times
    $t^4_{4k+3}$, $k\in\Z$.
\end{enumerate}
By Theorem \ref{dlmain}, this system solves a consensus problem.
The state trajectories of agents are shown in Fig.
\ref{switching}.

 \end{lz}
\section{Conclusion}
We presented an asynchronous consensus control strategy, which is
of obvious applications in realistic networks. By employing the
tools from the nonnegative matrix theory and graph theory, we
performed the convergence analysis of our consensus algorithm. The
introduction of communication topology $\Ge$ facilitated our
analysis and it established a connection between the actual
communication topology and our control strategy, and can be viewed
as the estimation of the actual topology. Examples were provided
to demonstrate the effectiveness of our theoretical results.

\appendix
\section{Lemmas}
\begin{yl}[\cite{J. Wolfowitz}, Lemma 1]\label{ylWolfowitz1}
If one ore more matrices in a product of stochastic matrices is
scrambling, so is the product.
\end{yl}
\begin{yl}[\cite{J. Wolfowitz}, Lemma 3]\label{ylWolfowitz3}
Let $A_1,A_2$ be stochastic matrices. If $A_2$ is an SIA matrix
and $A_1A_2\sim A_1$, then $A_1$ is a scrambling matrix.
\end{yl}
\begin{dy}
A nonnegative matrix $A\in \R^{n\times n}$ is said to be {\it
primitive} if it is irreducible and has only one eigenvalue of
maximum modulus.
\end{dy}

\begin{yl}[\cite{R.  Horn and C.  R.  Johnson}, pp.511, Corollary 8.4.8; pp.522, Problem
5]\label{yl848} Let $A\in \R^{n\times n}$ be nonnegative and
irreducible. If at least one main diagonal entry is positive, then
$A$ is primitive.
\end{yl}
\begin{yl}[\cite{R.  Horn and C. R.  Johnson}, pp.497, Lemma
8.2.7]\label{yl827} Let $A\in \R^{n\times n}$ be given, let
$\lambda\in \R$ be given, and suppose $\xi$ and $\zeta$ are
vectors such that
\begin{enumerate}
    \item $A\xi=\lambda \xi$;
    \item $A^T\zeta=\lambda \zeta$;
    \item $\xi^T\zeta=1$;
    \item $\lambda$ is an eigenvalue of $A$ with geometric multiplicity $1$;
    \item $|\lambda|=\rho(A)>0$; and
    \item $\lambda$ is the only eigenvalue of $A$ with modulus
    $\rho(A)$,
\end{enumerate}
where $\rho(A)$ is the spectral radius of $A$.
 Define $L=\xi \zeta^T$. Then
$(\lambda^{-1}A)^k=L+(\lambda^{-1}A-L)^k\to L$ as $k\to \infty$.
\end{yl}
\begin{yl}[\cite{W. Ren and R. W. Beard}, Lemma 3.4]\label{RWyl}
Let $A$ be a stochastic matrix. $\mathcal{G}(A)$ has a spanning
tree if and only if the eigenvalue $1$ of $A$ has algebraic
multiplicity equal to one.
\end{yl}

\begin{yl}[\cite{R.  Horn and C. R.  Johnson}, pp. 503, Theorem 8.3.1]\label{yl831}
If $A\in\R^{n\times n}$ and $A\geq 0$, then $\rho(A)$ is an
eigenvalue of $A$ and there is a nonnegative vector $f\geq 0$,
$f\not= 0$, such that $Af=\rho(A)f$.
\end{yl}

\end{document}